\documentclass[12pt]{article}%
\usepackage[applemac]{inputenc}
\usepackage{amsmath,amssymb}
\usepackage{amsfonts}
\usepackage[applemac]{inputenc}
\usepackage{amsmath,amssymb}
\usepackage{color}
\usepackage{color, colortbl}
\usepackage{amsmath}
\usepackage{amssymb}
\usepackage{graphicx}
\usepackage{tikz}
\usepackage{geometry}
\usetikzlibrary{arrows}
\usepackage{amsfonts}
\usepackage{amsmath,amssymb}
\usepackage[applemac]{inputenc}
\usepackage{color}
\usepackage{amsmath}
\usepackage{amssymb}
\usepackage{graphicx}
\usepackage{tikz}
\newtheorem{theorem}{Theorem}[section]
\newtheorem{theoremm}{Theorem}[subsection]
\newtheorem{lemma}[theorem]{Lemma}
\newcommand{\E}{E}
\newtheorem{proposition}[theorem]{Proposition}

\newtheorem{definition}[theorem]{Definition}

\newtheorem{example}[theorem]{Example}

\newtheorem{problem}[theorem]{Problem}

\newtheorem{definitionn}[theoremm]{Definition}

\newtheorem{examplee}[theoremm]{Example}

\newtheorem{remark}[theorem]{Remark}

\usetikzlibrary{arrows,shapes,automata,petri}
\newenvironment{myenumerate}{

\begin{enumerate}}{\end{enumerate}}
\newcommand{\dproof}{\noindent {Proof.} \quad}
\newcommand{\fproof}{\hfill $\square$ \bigskip}

\numberwithin{equation}{section}

\definecolor{LightCyan}{rgb}{0.88,1,1}

\def\RR{{\mathbb{ R}}}

\def\1B{\text{1\!\!I}}

\def\l{\langle}
\def\<{\langle}
\def\>{\rangle}
\def\P{\mathbb{P}}
\def\R{\mathbb{R}}
\def\l{\lambda}
\def\S{\mathcal{S}}
\def\E{\mathbb{E}}
\def\N{\mathbb{N}}
\def\L{\mathcal{L}}

\begin{document}

\title{The time-fractional heat equation driven by fractional time-space white noise}

\author{Rahma Yasmina Moulay Hachemi$^{1}$ \& Bernt \O ksendal$^{2}$
}
\date{19 February 2024
\vskip 0.5cm
\textcolor{blue}{In memory of Bent Fuglede}
}
\maketitle

\footnotetext[1]{%
Department of Mathematics, University of Saida, Algeria.\newline
Email: yasmin.moulayhachemi@yahoo.com}

\footnotetext[2]{%
Department of Mathematics, University of Oslo, Norway. 
Email: oksendal@math.uio.no.}


\begin{abstract}
We give an introduction to the  time-fractional stochastic heat equation driven by 1+d-parameter fractional time-space white noise, in the following two cases:
\begin{myenumerate}
    \item
With additive noise
  \item 
With multiplicative noise.
  \end{myenumerate}
  The fractional time derivative is interpreted as the Caputo derivative of order $\alpha \in (0,2)$ and we assume that the Hurst coefficient $H=(H_0,H_1,H_2, ...,H_d)$ of the time-space fractional white noise is in $(\tfrac{1}{2},1)^{1+d}$.

  We find an explicit expression for the unique solution in the sense of distribution of the equation in the additive noise case (i).
  
  In the multiplicative case (ii) we show that there is a unique solution in the Hida space $(\mathcal{S})^{*}$ of stochastic distributions and we show that the solution coincides with the solution of an associated fractional stochastic Volterra equation.Then we give an explicit expression for the solution of this Volterra equation.

A solution $Y(t,x)$ is called \emph{mild} if $E[Y^2(t,x)] < \infty$ for all $t,x$. For both the additive noise case and the multiplicative noise case we show that  if $\alpha \geq 1$ then the solution is mild if $d=1$ or $d=2$, while if $\alpha < 1$ the solution is not mild for any $d$.

The paper is partly a survey paper, explaining the concepts and methods behind the results. It is also partly a research paper, in the sense that some results appears to be new.

  \end{abstract}

  \section{Introduction}
\indent
The purpose of this paper is to study the time-fractional heat equation driven by multi-parameter fractional Brownian motion. We discuss both the additive noise case and the case with multiplicative noise. Using basic methods of Laplace and Fourier transforms combined with white noise calculus we obtain explicit solution formulas in the sense of distribution for any order $\alpha \in (0,2)$ of the time -derivative and for any space dimension $d$. Then we investigate for what values of $\alpha$ and $d$ the solution has finite second moment.

It was Niels Henrik Abel who first introduced the fractional derivative of a function.
In 1823 [1] he used the fractional derivative to solve the tautochrone (isochrone)
problem in mechanics.
Then in 1903 [22]  Gosta Magnus Mittag-
Leffler introduced what was later called the the Mittag-Leffler function $E_{\alpha}(z)$. Subsequently it became clear that this function has a connection to
the fractional derivative introduced by Abel.

As a modelling tool the fractional derivative of order $\alpha \in (0,2)$ is useful for many situations, e.g. in the
study of heat transfer and waves. For an example see Section 9. For many applications
of fractional derivatives we refer to the book by S. Holm [15].

 - The \emph{classical diffusion case} corresponds to $\alpha=1$. In that case our equation models the normal diffusion of
heat in a random or noisy medium, the noise being represented by the fractional
time-space white noise $W_H (t, x)$.

- The \emph{superdiffusion case} corresponds to $\alpha>1$. Then the equation models a situation where
the particles spread faster than in regular diffusion. This may occur for example in some
biological systems.

- The \emph{subdiffusion case}  occures when $\alpha <1$. Then the equation models a situation in which travel times of the
particles are longer than in the standard case. Such situations may occur for example in some transport systems. For more details see \cite{cherstry}.

Specifically, in this paper we study the following two  time-fractional stochastic heat equations driven by fractional time-space white noise:
\begin{align}
&\text{Additive noise:}\label{heat0}\nonumber\\
 	 &\frac{\partial^{\alpha}}{\partial t^{\alpha}}Y(t,x)=\lambda \Delta Y(t,x)+\sigma  W_{H}(t,x);\; (t,x)\in (0,\infty)\times \mathbb{R}^{d}.\\
&\text{Multiplicative noise:}\label{heat1}\nonumber\\
 	 &\frac{\partial^{\alpha}}{\partial t^{\alpha}}Y(t,x)=\lambda \Delta Y(t,x)+\sigma Y(t,x) \diamond W_{H}(t,x);\; (t,x)\in (0,\infty)\times \mathbb{R}^{d},
\end{align}
  	where $d\in\mathbb{N}=\{1,2,...\}$ and $\frac{\partial^{\alpha}}{\partial t^{\alpha}}$ is the Caputo derivative of order $\alpha \in (0,2)$, and $\l>0$ and $\sigma\in \mathbb{R}$ are given constants. 
Here 
\begin{equation}
  		\Delta Y =\sum_{j=1}^{d}\frac{\partial^{2}Y}{\partial x_{j}^{2}}(t,x)
  	\end{equation} 
denotes the Laplacian operator acting on $x$, $W_{H}(t,x)$ is multiparameter fractional  white noise, defined by
\begin{equation}
  	W_{H}(t,x)=\frac{\partial}{\partial t}\frac{\partial^{d}B_{H}(t,x)}{\partial x_{1}...\partial x_{d}},
  \end{equation}

\noindent where $B_{H}(t,x)=B(t,x,\omega); t\geq 0, x \in \R^d, \omega \in \Omega$ is time-space fractional Brownian motion with Hurst parameter $H=(H_0, H_1, ... , H_d) \in (\tfrac{1}{2},1)^{d+1}$ with probability law $\P$, and $\diamond$ denotes the Wick product.

The boundary conditions are
\begin{align}
    Y(0,x)&=\delta_0(x)\text{ (the point mass at  } 0), \label{1.4}\\
    \lim_{x \rightarrow +/- \text{ }\infty}Y(t,x)&=0.\label{1.5}
\end{align}
 We consider the equation \eqref{heat0} in the sense of distribution, and we find an explicit expression for the $\mathcal{S}'$-valued solution $Y(t,x)$, where $\mathcal{S}'$ is the space of tempered distributions. See Theorem \ref{add}.\\
 
The equation \eqref{heat1} is considered in the Hida space $(\mathcal{S})^{*}$ of stochastic distributions. We show that $Y(t,x)$ solves \eqref{heat1} if and only if $Y(t,x)$ solves a Volterra equation driven by the same fractional Brownian motion. See Theorem \ref{mult1}.

Then we find an explicit expression for the solution of such Volterra equations. See Theorem \ref{vol3}.

The solution is called  \emph{mild} if $Y(t,x) \in L^2(\P)$ for all $t,x$.  It is well-known (see e.g. Y. Hu \cite{Hu}) that in the classical case with $\alpha = 1$, the solution is mild if and only if the space dimension $d=1$. We prove that if $\alpha \in (1,2)$ the solution is mild if $d=1$ or $d=2$.
 If $\alpha < 1$ we prove that the solution is not mild for any $d$. This holds for both the additive noise and the multiplicative noise case and for all Hurst vectors $H=(H_0,H_1, ... ,H_d) \in (\tfrac{1}{2},1)^{1+d}$. See Theorem \ref{mild1} and Theorem \ref{mild2}. 
 
 A crucial base for our results is our $L^2(P)$-estimate for WIS-stochastic integrals with respect to fractional Brownian motion (Theorem \ref{th6.1}).This result appears to be new.\\

There are many papers dealing with various forms of stochastic fractional differential equations. Some papers which are related to ours are:  \\
 -  In the paper by Kochubel et al \cite{KKd} the fractional heat equation corresponding to random time change in Brownian motion is studied.\\
 - The papers by Bock et al \cite{BGO}, \cite{Bd} 
 are considering stochastic equations driven by grey Brownian motion. \\
  - The paper Roeckner et al \cite{LRd} proves existence and uniqueness of general time-fractional linear evolution equations in the Gelfand triple setting.\\
  -The paper Yalin et al \cite{Y} study the time-regularity of the paths of solutions to stochastic
  partial differential equations  driven by additive
  infinite-dimensional fractional Brownian noise.\\
  -The paper by Binh et al  \cite{Bi} study the spatially-temporally H\o lder continuity of mild random field solution of space time fractional stochastic heat equation driven by colored noise.\\
  - The paper  which is closest to our paper is Chen et al \cite{CHN}, where 
a comprehensive discussion is given of a general fractional stochastic heat equations with multiplicative noise, and with fractional derivatives in both time and space, is given. In that paper the authors prove existence and uniqueness results as well as regularity results of the solution, and they give sufficient conditions on the coefficients and the space dimension $d$, for the solution to be a random field. Their method is based on the use of fundamental solutions of the corresponding deterministic equations.

In comparison, our paper  is dealing with both the additive and multiplicative noise case. We consider a more special class of fractional heat equations than in Chen et al \cite{CHN} and this allows us to use basic methods of Laplace and Fourier transform and to find explicit solution formulas in the sense of distributions.
In Section 5 we prove a new $L^2(P)$-estimate of independent interest for fBm integrals. Then we use this estimate combined with the explicit solution formulas to find conditions under which the solutions are random fields in $L^2(\P)$. Our paper may be regarded as an extension of the paper \cite{MO}, where only the case with additive classical Brownian sheet noise (not fractional noise) is considered.

We refer to Holm \cite{H}, Ibe \cite{I}, Kilbas et al \cite{Kilbas} and Samko et al \cite{Samko} for more information about fractional calculus and their applications.

\section{Preliminaries}
In this section we summarise some of the background notation and results needed in this paper.
\subsection{The space of tempered distributions}
For the convenience of the reader we recall some of the basic properties of the Schwartz space $\mathcal{S}$ of rapidly decreasing smooth functions and its dual, the space $\mathcal{S}'$ of tempered distributions.
Let $n$ be a given natural number. Let $\mathcal{S}=\mathcal{S}(\mathbb{R}^n)$\label{simb-028} be the
space of rapidly decreasing smooth real
functions $f$
on $\mathbb{R}^n$
equipped with the family of seminorms:\label{simb-029} 
\begin{equation*}
\Vert f \Vert_{k,\alpha} := \sup_{y \in \mathbb{R}^n}\big\{ (1+|y|^k) \vert
\partial^\alpha f(y)\vert \big\}< \infty,
\end{equation*}
where $k = 0,1,...$, $\alpha=(\alpha_1,...,\alpha_n)$ is a multi-index with $%
\alpha_j= 0,1,...$ $(j=1,...,n)$ and\label{simb-030} 
\begin{equation*}
\partial^\alpha f := \frac{\partial^{|\alpha|}}{\partial
y_1^{\alpha_1}\cdots \partial y_n^{\alpha_n}}f
\end{equation*}
for $|\alpha|=\alpha_1+ ... +\alpha_n$.

Then
$\mathcal{S}=\mathcal{S}(\mathbb{R}^n)$ is a
Fr\'echet space.

Let $\mathcal{S}^{\prime }=\mathcal{S}^{\prime }(\mathbb{R}^{n})$\label%
{simb-031} be its dual, called the space of \emph{tempered distributions}. 
\index{tempered distributions} Let $\mathcal{B}$ denote the family of all
Borel subsets of $\mathcal{S}^{\prime }(\mathbb{R}^{n})$ equipped with the
weak* topology. If $\Phi \in \mathcal{S}^{\prime }$ and $f \in \mathcal{%
S}$ we let \label{simb-033} 
\begin{equation}
\Phi (f) \text{ or } \langle \Phi ,f \rangle  \label{3.1}
\end{equation}%
denote the action of $\Phi$ on $f$.

\begin{example}
\begin{itemize}
\item
{(Evaluations)}
For $y \in \R$ define the function $\delta_y$ on $\S(\R)$ by $\delta_y(\phi)=\phi(y)$. Then $\delta_y$ is a tempered distribution.\\
\item
{(Derivatives)} Consider the function $D$, defined for $\phi \in \S(\mathbb{R})$ by $D[\phi]=\phi^{\prime}(y)$. Then  $D$ is a tempered distribution. \\
\item
{(Distributional derivative)}\\
 Let $T$ be a tempered distribution, i.e. $T \in \S^{'}(\mathbb{R}) $. We define the distributional derivative $T^{'}$ of $T$ by
 $$ T^{'}[\phi]=-T[\phi^{'}]; \quad \phi \in \S.$$
 Then $T^{'}$ is again a tempered distribution.
 \end{itemize}
 \end{example}
 
 In the following we will apply this to the case when $n=1+d$ and $y=(t,x) \in \R \times \R^d$.
 \subsection{The Mittag-Leffler functions}
 \begin{definition}(The two-parameter Mittag-Leffler function)
 The Mittag-Leffler function of two parameters $\alpha,\; \beta$ is denoted by $E_{\alpha,\beta}(z)$ and defined by:
\begin{equation}
    E_{\alpha,\beta}(z)=\sum_{k=0}^{\infty}\frac{z^{k}}{\Gamma(\alpha k+\beta)}
\end{equation}
where $z,\; \alpha,\; \beta\in \mathbb{C}$ (the set of complex numbers),\; $Re(\alpha)>0\; and\; Re(\beta)>0,$ and $\Gamma$ is the Gamma function.\\

 For $\beta=1$ we obtain
the Mittag-Leffler function of one parameter $\alpha$ denoted by $E_{\alpha}(z)$. This function is defined as
\begin{equation}
    E_{\alpha}(z)=\sum_{k=0}^{\infty}\frac{z^{k}}{\Gamma(\alpha k+1)}
\end{equation}
where $z,\; \alpha\in \mathbb{C},\; Re(\alpha)>0.$
\end{definition}

\begin{remark} Note that $E_{\alpha}(z)= E_{\alpha,1}(z)$ and that 
\begin{align}
E_{1}(z)=\sum_{k=0}^{\infty}\frac{z^{k}}{\Gamma( k+1)}=\sum_{k=0}^{\infty}\frac{z^{k}}{k!} =e^{z}.
\end{align}
\end{remark}

\subsection{The (Abel-)Caputo fractional derivative}
  In this section we present the definitions and some properties of the Caputo derivatives.
  \begin{definition}
  The (Abel-)Caputo fractional derivative of order $\alpha > 0$ of a function $f$ such that $f(x)=0$ when $x<0$ is denoted by  $D^{\alpha} f (x)$ or $\frac{d^{\alpha}}{dx^{\alpha}} f(x)$ and defined  by 
  \begin{align}\label{caputo1}
  D^{\alpha}f(x):& =
  \begin{cases}
  \frac{1}{\Gamma(n-\alpha)}\int_0^x \frac{f^{(n)}(u)du}{(x-u)^{\alpha +1 -n}}; \quad n-1 < \alpha < n\\
  \frac{d^n}{dx^n}f(x); \quad \alpha =n.
  \end{cases}
  \end{align}
  Here $n$ is an smallest integer greater than or equal to $\alpha$.\\
  
  \noindent If $f$ is not smooth these derivatives are interpreted in the sense of distributions.
  \end{definition}
  
  \begin{example}
  If $f(x)=x$ and $\alpha \in (0,1)$ then
 
  \begin{align}
 D^{\alpha}f(x)=\frac{ x^{1-\alpha} }{(1-\alpha)\Gamma(1-\alpha)}.
 \end{align} 
 In particular, choosing $\alpha=\tfrac{1}{2}$ we get
 \begin{align}
 D^{\tfrac{1}{2}}f(x)=\frac{2 \sqrt{x}}{\sqrt{\pi}}.
 \end{align} 
  \end{example}

\subsection{Laplace transform of Caputo derivatives} 
\begin{itemize}
\item
The 1-dimensional case.\\
Recall that the (1-dimensional) Laplace transform $L$ is defined by 
 
\begin{equation}
	Lf(s)=\int_{0}^{\infty}e^{-st}f(t)dt=\widetilde{f}(s)
	\end{equation}
	for all $f$ such that the integral converges.

 Some of the properties of the Laplace transform that we will need are:

  \begin{align}
    &L[ \frac{\partial ^{\alpha}}{\partial t^{\alpha}}f(t)](s)=s^{\alpha}(L f)(s)-s^{\alpha-1}f(0) \label{L1}\\       
     &L[E_{\alpha}(bx^{\alpha})](s)  = \frac{s^{\alpha -1}}{s^{\alpha}-b}\label{L2}\\
     &L[x^{\alpha-1}E_{\alpha,\alpha}(-b x^{\alpha})](s)=\frac{1}{s^{\alpha}+b}\label{L3}
     \end{align}
     Recall that the convolution $f\ast g$ of two functions $f,g: [0,\infty) \mapsto \mathbb{R}$  is defined by
\begin{align}
(f \ast g)(t)=\int_0^t f(t-r)g(r) dr; \quad t \geq 0.
\end{align}

The convolution rule for Laplace transform states that $$L\left( \int_{0}^{t}f(t-r)g(r)dr\right) (s)=Lf(s)Lg(s),$$ 
or 
\begin{equation}\label{12aa}
	\int_{0}^{t}f(t-r)g(r)dr=L^{-1}\left( Lf(s)Lg(s)\right) (t).
\end{equation}

\item
{The multidimensional case}\\
If $g: \mathbb{R}_{+}^d \mapsto \mathbb{R}$ is bounded, we define the \emph{multidimensional }Laplace transform $\L$ of $g$ by

\begin{equation}
	\L\{g\}(y)=\int_{\mathbb{R}^d_{+}}e^{-xy}g(x)dx=:\widehat{g}(y);\quad y\in \mathbb{R}_{+}^d; 
	\end{equation}
 where $xy=x_1y_2 + x_2y_2 + ... +x_dy_d$ when $x=(x_1,x_2, ... ,x_d), y=(y_1,y_2, ... ,y_d).$\\
 In this case the convolution $f\ast g$ of two functions $f,g: [0,\infty)^d \mapsto \mathbb{R}$  is defined by
\begin{align}
(f \ast g)(x)=\int_0^x f(x-z)g(z) dz; \quad x \in \mathbb{R}_{+}^d,
\end{align}
where we use the notation
$\int_0^x = \int_0^{x_1} \int_0^{x_2} ... \int_0^{x_d}$ if $x=(x_1,x_2, ... x_d).$\\

The convolution rule for the multidimensional Laplace transform states that $$\mathcal{L}\left( \int_{0}^{x}f(x-z)g(z)dz\right) (y)=\L f(y) \L g(y),$$ 
or 
\begin{equation}\label{12A}
	\int_{0}^{x}f(x-z)g(z)dz=\L ^{-1}\left( \L f(y)Lg(y)\right) (x).
\end{equation}

\end{itemize}

 \section{Time-space white noise}
In this section we briefly review the basic notation and results of time-space white noise. 

Let $n$ be a fixed natural number. Later we will set $n= 1 + d$.
Define $\Omega={\cal S}'(\RR^n)$, equipped
with the weak-star topology. This space will be the base of our basic
probability space, which we explain in the following:
\vskip 0.3cm
As events we will use the family $\mathcal{F}={\cal
B}({\cal S}'(\RR^n))$ of Borel subsets of ${\cal S}'(\RR^d)$, and our
probability measure $\P$ is defined by the following result:

\begin{theorem}{\bf (The Bochner--Minlos theorem)}\\
There exists a unique probability measure $\P$ on ${\cal B}({\cal S}'(\RR^n))$
with the following property:
$$\E[e^{i\langle\cdot,\phi\rangle}]:=\int\limits_{\cal S'}e^{i\langle\omega,
\phi\rangle}d\mu(\omega)=e^{-\tfrac{1}{2} \Vert\phi\Vert^2};\quad i=\sqrt{-1}$$
for all $\phi\in{\cal S}(\RR^n)$, where
$\Vert\phi\Vert^2=\Vert\phi\Vert^2_{L^2(\RR^n)},\quad\langle\omega,\phi\rangle=
\omega(\phi)$ is the action of $\omega\in{\cal S}'(\RR^n)$ on
$\phi\in{\cal S}(\RR^n)$ and $\E=\E_{\P}$ denotes the expectation
with respect to  $\P$.
\end{theorem} 
We will call the triplet $({\cal S}'(\RR^n),{\cal
B}({\cal S}'(\RR^n)),\P)$ the {\it  white noise probability
space\/}, and $\P$ is called the {\it white noise probability measure}.

The measure $\P$ is also often called the (normalised) {\it Gaussian
measure\/} on ${\cal S}'(\RR^n)$. It is not difficult to prove that if $\phi\in L^2(\RR^n)$ and
we choose $\phi_k\in{\cal S}(\RR^n)$ such that $\phi_k\to\phi$ in $L^2(\RR^n)$,
then
$$\langle\omega,\phi\rangle:=\lim\limits_{k\to\infty}\langle\omega,\phi_k\rangle
\quad\text{exists in}\quad L^2(\P)$$
and is independent of the choice of $\{\phi_k\}$. In
particular, if we define
$$\widetilde{B}(x):=\widetilde{B}(x_1,\cdots,x_n,\omega)=\langle\omega,\chi_
{[0,x_1]\times\cdots\times[0,x_n]}\rangle; \quad  
x=(x_1,\cdots,x_n)\in\RR^n,$$
where $[0,x_i]$ is interpreted as $[x_i,0]$ if $x_i<0$,
then $\widetilde{B}(x,\omega)$ has an $x$-continuous version $B(x,\omega)$, which
becomes an \emph{$n$-parameter Brownian motion}, in the following sense:

By an \emph{$n$-parameter Brownian motion} we mean a family
$\{B(x,\cdot)\}_{x\in\RR^n}$ of random variables on a probability space
$(\Omega,{\cal F},\P)$ such that
\begin{itemize} 
\item
$B(0,\cdot)=0\quad\text{almost surely with respect to } \P,$
\item
$\{B(x,\omega)\}$ is a continuous and Gaussian stochastic process 
\item
For all $x=(x_1,\cdots,x_n)$,
$y=(y_1,\cdots,y_n)\in\RR_+^n$,
$B(x,\cdot),\,B(y,\cdot)$ have the covariance
$\prod_{i=1}^n x_i\wedge y_i$. 
\end{itemize} 

It can be proved that the process $\widetilde{B}(x,\omega)$ defined above has a modification $B(x,\omega)$ which satisfies all these properties.
This process $B(x,\omega)$ then becomes an \emph{$n$-parameter Brownian motion}. 
\vskip0.3cm
We remark that for $n=1$ we get the classical (1-parameter) Brownian motion
$B(t)$ if we restrict ourselves to $t\geq 0$. For $n \geq 2$ we get what is often
called the \emph{Brownian sheet}.
\vskip0.3cm
With this definition of Brownian motion it is natural to define the
$n$-parameter Wiener--It\^{o} integral of $\phi\in L^2(\RR^n)$ by
$$\int\limits_{\RR^n}\phi(x)dB(x,\omega):=\langle\omega,\phi\rangle;\quad
\omega\in{\cal S}'(\RR^n).$$
We see that by using the Bochner--Minlos theorem we have obtained an easy construction of
$n$-parameter Brownian motion that works for any parameter dimension $n$. Moreover, we get a representation of the space $\Omega$ as the Fr\' echet space $\mathcal{S}'(\R^n)$. This is an advantage in many situations, for example in the construction of the Hida-Malliavin derivative, which can be regarded as a stochastic gradient on $\Omega$.
See e.g. \cite{Dop} and the references therein.

In the following we put $n=1+d$ and let 
 $$B(t,x)=B(t,x,\omega); t \geq 0, x \in \R^d, \omega \in \Omega$$ 
 denote the (1-dimensional) time-space Brownian motion (also called the Brownian sheet) with probability law $\P$. Since this process is $(t,x)$-continuous a.s., we can for a.a. $\omega \in \Omega$ define its derivatives with respect to $t$ and $x$ in the sense of distributions.
 Thus we define the time-space white noise $W(t,x)=W(t,x,\omega)$ by
 
\begin{equation}
  	W(t,x)=\frac{\partial}{\partial t}\frac{\partial^{d}B(t,x)}{\partial x_{1}...\partial x_{d}}. \label{WN}
  \end{equation}
In particular, for $d=1$ and $x_1=t$ and get the familiar identity
$$W(t)=\frac{d}{dt}B(t)\hbox{ in } \S^{'}.$$ 

The definition \eqref{WN} can also be interpreted as an element of the Hida space $(\S)^*$ of \emph{stochastic distributions}, and in that setting it has been proved (see Lindstr\o m, \O . \& Ub\o e \cite{LOU} and Benth \cite{B}) that the Ito-Skorohod integral with respect to $B(dt,dx)$ can be expressed as 
\begin{align}
\int_0^T \int_{\mathbb{R}^d} f(t,x,\omega) B(dt, dx)=\int_0^T \int_{\mathbb{R}^d} f(t,x,\omega) \diamond W(t, x) dt dx,
\end{align}
where $\diamond$ denotes the Wick product.\\
\emph{This is the interpretation we are using in \eqref{heat1} and elsewhere in this paper.}

In particular, if $f(t, x,\omega)=f(t, x)$ is deterministic, we can drop the Wick product and get
\begin{align}
\int_0^T \int_{\mathbb{R}^d} f(t, x) B(dt, dx)=\int_0^T \int_{\mathbb{R}^d} f(t,x)  W(t,x) dt dx.
\end{align}

\section{Fractional Brownian motion and fractional white noise calculus}

There are several ways to introduce fractional Brownian motion and fractional white noise calculus. One commonly used approach due to Hu et al \cite{Hu} is recalled in the Appendix. For our purpose it is more convenient to follow the approach presented by Elliott and van der Hoek (see \cite{EvdH} and the references therein, including \cite{L}). Hence we introduce the multi-parameter (time-space) fractional Brownian motion  as obtained by applying an operator $M$ to the classical multi-parameter Brownian motion.  This approach has several advantages. For example, all the different fBm's corresponding to different Hurst coefficients $H \in (0,1)$ are defined on the same white noise probability space $(\Omega, \mathbb{F},P)$ introduced in Section 3. There are also several computational advantages, as shown later. See also e.g. \cite{DO}.\\

In the following we let $f: \mathbb{R}^n \mapsto \mathbb{R}$ be a function such that its Fourier transform
\begin{equation}
\mathcal{F}(f)(y):= \widehat{f}(y):= (2 \pi)^{-\tfrac{n}{2}} \int_{\mathbb{R}^n} e^{-ixy}f(x)dx; \quad y\in \mathbb{R}^n,\quad ( i=\sqrt{-1})
 \end{equation}
exists.
Let $\mathcal{F}^{-1}$ denote the inverse Fourier transform, defined by
 \begin{align}
 \mathcal{F}^{-1} g(x):= (2 \pi)^{-\tfrac{n}{2}}\int_{\mathbb{R}^n} e^{ixy} g(y) dy,
 \end{align}
 for all functions $g$ such that the integral converges.

\begin{definition}(\cite{EvdH},(A10))
Let $H=(H_1, H_2, ... ,H_n) \in (\tfrac{1}{2},1)^n $ and let $\mathcal{S}$ be the Schwartz space of rapidly decreasing smooth functions on $\mathbb{R}^n$. The operator $M=M^{(H)}: \mathcal{S} \mapsto \mathcal{S}$ is defined by
\begin{align}
Mf(x)&=\mathcal{F}^{-1}(\prod_{j=1}^n |y_j |^{\frac{1}{2}-H_j} (\mathcal{F}f)(y))(x)\\
&=(2 \pi)^{-\tfrac{d}{2}}\int_{\mathbb{R}^n} e^{ixy} (\prod_{j=1}^n |y_j|^{\frac{1}{2} - H_j}(2 \pi)^{-\tfrac{d}{2}} \int_{\mathbb{R}^n} e^{-ixy}f(x)dx) dy
\end{align}
If $H_j > \frac{1}{2} $ for all $j$ then the operator $M$ can also be defined as follows:
\begin{align}
M_Hf(x)=  \int_{\mathbb{R}^n} \prod_{j=1}^n C_{H_j}| y_j | ^{H_j-\frac{3}{2}} f(x+y) dy; \quad f \in \mathcal{S},\label{MHa}
\end{align}
where
\begin{equation}
C_{H_j}=[2\Gamma(H_j-\tfrac{1}{2}) \cos(\tfrac{\pi}{2}(H_j-\tfrac{1}{2}))] ^{-1},
\end{equation}
$\Gamma$ being the classical Gamma function.
\end{definition}
\emph{From now on we assume that $H \in (\frac{1}{2},1)^n$.}\\

Since $H \in (\frac{1}{2},1)^n$ is fixed, we will suppress the index $H$ and write $M^{(H)} = M$ in the following.\\
We can in a natural way extend the operator $M$  from $\mathcal{S}$ to the space
\begin{align}
L_{H}^2(\mathbb{R}^n):= \{ f: \mathbb{R}^n \mapsto \mathbb{R}; Mf \in L^2(\mathbb{R}^n )\}.
\end{align}
Note that $L_{H}^2 (\mathbb{R}^n)$ is a Hilbert space when equipped with the inner product
\begin{equation}
(f,g)_H := (f,g)_{L_H^2(\mathbb{R}^n)}=(Mf,Mg)_{L^2 (\mathbb{R}^n)};\quad f,g \in L^2_{H}(\mathbb{R}^n).
\end{equation}

The following result is useful:
\begin{lemma}
Let $\varphi \in L^2(\mathbb{R}^n)$ and $\psi \in L^2_{H}(\mathbb{R}^n)$. Then
\begin{align}
(\varphi,M\psi)_{L^2(\mathbb{R}^n)}=(M\varphi,\psi)_{L^2(\mathbb{R}^n)}
\end{align}

\end{lemma}

\dproof
Since
\begin{equation}
Mf(x)= \int_{\mathbb{R}^n} \prod_{j=1}^n C_{H_j}| y_j | ^{H_j-\frac{3}{2}} f(x+y) dy; \quad f \in \mathcal{S},\label{MH}
\end{equation}
we get by the Fubini-Tonelli theorem
\begin{align*}
(\varphi, M\psi)_{L^2(\mathbb{R}^n)}&= \int_{\mathbb{R}^n} \varphi(x)M\psi(x)dx= \int_{\mathbb{R}^n}\varphi(x) \int_{\mathbb{R}^n} \prod_{j=1}^n C_{H_j}| y_j | ^{H_j-\frac{3}{2}} \psi(x+y) dy dx\\
&= \int_{\mathbb{R}^n}\int_{\mathbb{R}^n}\varphi(z-y)\prod_{j=1}^n C_{H_j}| y_j | ^{H_j-\frac{3}{2}}\psi(z) dy dz\\
&=  \int_{\mathbb{R}^n}\int_{\mathbb{R}^n}\psi(z) \prod_{j=1}^n C_{H_j}| y_j | ^{H_j-\frac{3}{2}} \varphi(z+y) dy dz=(M\varphi,\psi)_{L^2(\mathbb{R})},
\end{align*}
which completes the proof.

\fproof

The following property of the operator $M$ will be used in the sequel.

\begin{lemma}\label{l1} 
The operator $M$ satisfies following property:
\begin{align} \label{eq3.5a}
M^2f(x)= \int_{\mathbb{R}^n} \int_{\mathbb{R}^n} \prod_{j=1}^n C_{H_j}^2| y_j z_j| ^{H_j-\frac{3}{2}}  f(y+z+x)dy dz.
\end{align}
\end{lemma}

\dproof 
Since
$$Mf(x)=\int_{\mathbb{R}^n}\prod_{j=1}^n C_{H_j}| y_j | ^{H_j-\frac{3}{2}} f(x+y)dy$$
we get
\begin{align*}
M^2f(x)&=M(Mf)(x)= \int_{\mathbb{R}^n}\prod_{j=1}^n C_{H_j}| z_j | ^{H_j-\frac{3}{2}}  Mf(x+z) dz \\
&= \int_{\mathbb{R}^n} \prod_{j=1}^n C_{H_j}| z_j | ^{H_j-\frac{3}{2}} \int_{\mathbb{R}^n}\prod_{j=1}^n C_{H_j}| y_j | ^{H_j-\frac{3}{2}}f(x+z+y)dy dz\\
&=\int_{\mathbb{R}^n} \int_{\mathbb{R}^n}\prod_{j=1}^n C_{H_j}^2| y_j z_j | ^{H_j-\frac{3}{2}} f(y+z+x)dy dz.
\end{align*}
This completes the proof.
\fproof

Results from previous lemmas can be gathered in the extended equality which is given in following lemma.

\begin{lemma} \label{5.4} Let  $f,g:\mathbb{R}^n\rightarrow \mathbb{R}$. Then
	\begin{align} \label{eq l21}
	\int_{\mathbb{R}^n}f(x)	M^2 g(x)dx= \int_{\mathbb{R}^n}g(x)M^2 g(x)dx=\int_{\mathbb{R}^n} Mf(x) Mg(x) dx.
	\end{align}
\end{lemma}

\dproof 
From Lemma \ref{l1} it follows that
\begin{align}
\int_{\mathbb{R}^n} f(x)	M^2 g(x)dx&=\int_{\mathbb{R}^n}f(x) \int_{\mathbb{R}^n} \int_{\mathbb{R}^n} \prod_{j=1}^n C_{H_j}^2| y_j z_j| ^{H_j-\frac{3}{2}}  g(y+z+x)dy dz dx\label{eq3.5}
\end{align}.

If we substitute

\begin{equation*} 
	\left\{ 
        \begin{array}{l} 
            u=x+y+z\\
            v=y\\
            w=z\\
        \end{array}
        \right. 	\Leftrightarrow  \begin{array}{l} 
            x=u-v-w=: X(u,v,w)\\
            y=v=:Y(u,v,w)\\
            z=w=:Z(u,v,w),\\
        \end{array}
    \end{equation*}    
    then the Jacobian is 
    \begin{equation}
    \begin{vmatrix} 
\frac{\partial X}{\partial u} & \frac{\partial X}{\partial v} & \frac{\partial X}{\partial w} \\ 
\frac{\partial Y}{\partial u} & \frac{\partial Y}{\partial v} & \frac{\partial Y}{\partial w} \\ 
\frac{\partial Z}{\partial u} & \frac{\partial Z}{\partial v} & \frac{\partial Z}{\partial w} 
\end{vmatrix} =   
    \begin{vmatrix}
1 & -1 & -1 \\ 
0 & 1 & 0\\ 
0 & 0 & 1
\end{vmatrix}=1,
\end{equation}	

and we obtain
\begin{align}
&\int_{\mathbb{R}^n}f(x)	M^2g(x) dx=\int_{\mathbb{R}^n}\int_{\mathbb{R}^n}\int_{\mathbb{R}^n}\prod_{j=1}^n C_{H_j}^2| u_j w_j| ^{H_j-\frac{3}{2}}  f(u-v-w)g(u)du dv dw\nonumber\\
&=\int_{\mathbb{R}^n}\int_{\mathbb{R}^n}\int_{\mathbb{R}^n} g(u)\prod_{j=1}^n C_{H_j}^2| u_j w_j| ^{H_j-\frac{3}{2}}  f(u-v-w)g(u)du dv dw\nonumber\\	&=\int_{\mathbb{R}^n}g(u)	M^2 f(u)du.
\end{align}
This proves the first identity. A similar argument proves the second identity.
\fproof

\medskip

Now we can give a formal definition of multi-parameter fractional Brownian motion:

\begin{definition}{(Fractional Brownian motion revisited)}\\
For $x \in \mathbb{R}^n$ define
\begin{equation} \label{eq4.5a}
 B_{H}(x):= B_{H}(x,\omega):=\langle \omega, M_x(\cdot) \rangle ; \quad \omega \in \Omega 
\end{equation}
where we have put 
\begin{align}
M_x (y):=M(\chi_{[0,x_1]\times [0,x_2] \times ... \times [0,x_n]})(y)=\prod_{j=1}^n C_{H_j}^{-1}\Big(\frac{x_j-y_j}{|x_j-y_j|^{\frac{3}{2}-H}}+\frac{y_j}{|y_j|^{\frac{3}{2}-H}}\Big);\quad y,x \in \mathbb{R}^n.
\end{align}

Then $B_{H}$ is a Gaussian process, and by the properties of the classical Bm we have 
\begin{equation}
E[B_{H} (x)]=B_{H}(0)=0,
\end{equation}
and (see \cite{EvdH},(A10))
\begin{equation}
E[B_{H}(x) B_{H}(y)]=  \prod_{j=1}^{d}C_{H_j}\big( |x_{j}|^{2H_{j}}+|y_{j}|^{2H_{j}}-|x_{j}-y_{j}|^{2H_{j}}\big)]; \quad x,y \in \mathbb{R}^n.
\end{equation}

By the Kolmogorov continuity theorem the process $B_{H}$ has a continuous version, which we will also denote by $B_{H}.$ Except for a constant, this process coincides with the fractional Brownian motion $B_H(x)$ introduced in the Appendix.
From now on we call this process 
 \emph{ the fractional Brownian motion (fBm) with Hurst parameter $H$} and denote it by $B_H(x)$.
\end{definition}

\subsection{Fractional white noise calculus}
By \eqref{eq4.5a} it follows by approximating $f \in \mathcal{S}$ by step functions that
\begin{equation}\label{int}
 \int_{\mathbb{R}^n} f(x) dB_{H}(x)=\langle \omega,Mf \rangle =\int_{\mathbb{R}^n} Mf(x) dB(x)
\end{equation}
and hence
\begin{equation}
B_{H}(x)=\int_{\mathbb{R}^n} \prod_{j=1}^n \chi_{[0,x_j]}(y_j) dB_{H}(y)= \int_{\mathbb{R}^n}M_{x}(y)dB(y)=\langle \omega, M_x \rangle.
\end{equation}
We know that the action of $\omega \in \mathcal{S}'(\mathbb{R}^n)$ extends from $\mathcal{S}(\mathbb{R}^n)$ to $L^2(\mathbb{R}^n),$ and that if $f \in \mathcal{S}$ then $Mf \in L^2(\mathbb{R}^n).$ Therefore we can define $M$ on $\mathcal{S}'(\mathbb{R}^n)$ by setting
\begin{equation}
\langle M\omega ,\varphi \rangle= \langle \omega, M\varphi \rangle; \varphi \in \mathcal{S}(\mathbb{R}^n), \omega \in \mathcal{S}'(\mathbb{R}^n).
\end{equation}
Using this, we can write
\begin{align}
B_{H}(x,\omega)=\langle \omega, M_x \rangle= \langle M\omega, \prod_{j=1}^n\chi_{[0.x_j](\cdot)} \rangle =B(x,M\omega).
\end{align}

We define the \emph{fractional white noise} 
$W_{H}(x)=W_{H}(x,\omega)$ by
\begin{align}
W_{H}(x,\omega)= \frac{\partial^{n}}{\partial x_{1}...\partial x_{n}}B_{H}(x) \text{ (derivative in } (\mathcal{S})^{*} ).
\end{align}
By the above it follows that 

\begin{align}
W_{H}(x,\omega)=W(x,M\omega),
\end{align}
 where 
\begin{align}
W(x)=W(x,\omega)= \frac{\partial^{n}}{\partial x_{1}...\partial x_{n}}B(x)  \text{ is the classical white noise in } (\mathcal{S})^{*}.
\end{align}
Let $L_H^2(\R^n)$ denote the set of measurable processes $f(x)=f(x,\omega): \R^n \times \Omega \mapsto \R$ such that 
\begin{align}
E[\int_{\R^n} (M f)^2(x)dx] < \infty.
\end{align}
Then we define the \emph{Wick-Ito-Skorohod} (WIS)-integral  of $f \in L_H^2(\R^n)$  as follows:

\begin{definition}  \label{wis}
(The fractional Wick--It\^o--Skorohod (WIS) integral) \\
Let  $f : \mathbb{R}^n \rightarrow (\mathcal{S})^{*}$ be such that $f (x)  \diamond W_{H}(x) $ is $dx$-integrable in $ (\mathcal{S})^{*}$. Then we say that  $f$ is
$dB_{H}$-integrable and we define the Wick--It\^o--Skorohod (WIS) integral of $f (\cdot) = f(\cdot, \omega)$ with respect to $B_{H}$ by
\begin{align}
\int_{\mathbb{R}^n} f(x) dB_{H}(x)= \int_{\mathbb{R}^n} f(x) \diamond W_{H}(x) dx,
\end{align}
where the integral on the right is the (Pettis) Lebesgue integral with values in $(\mathcal{S})^{*}.$
\end{definition}

By \eqref{int} and Lemma \ref{5.4} we obtain the following isometry:
\begin{theorem} \label{iso1}
\begin{align}
    E[(\int_{\mathbb{R}^n} f(x) dB_{H}(x))^2]= E[\int_{\mathbb{R}^n} (Mf)^2(x) dx]= E[\int_{\mathbb{R}^n} f(x) M^2f(x) dx]
    \end{align}
\end{theorem}

\section{An $L^2(P)$-estimate for fBm -integrals}	
	The following theorem gives a fundamental  $L^2(P)$-inequality for WIS-integrals with respect to fBm.
	It is a multi-parameter extension of Theorem 3.1 in \cite{DO}.
	\begin{theorem} \label{th6.1}
	Let $f(x)=f(x,\omega)$   be a measurable  stochastic process such that
	
	$$E\left[\int_{\R^n} (Mf)^2(x)dx\right]<\infty.$$
	Assume that $f$ is supported on the set 
 $S=[0,\bar{x}_1] \times [0,\bar{x}_2] \times ... \times [0,\bar{x}_n]$. 
	Define
	
	\begin{align}
	Y&=\int_{0}^{\bar{x}} f(x) dB_{H}(x)=\int_{S} f(x) dB_{H}(x)\\
 &=\int_0^{\bar{x}_1} \int_0^{\bar{x}_2} ... \int_0^{\bar{x}_n} f(x_1,x_2, ... , x_n)dB_{H}(x_1, x_2, ... , x_n),
	\end{align}
	where the integral is interpreted in the WIS sense.
	
	Then
	\begin{equation}
	E[Y^2]\leq (\int_{S}E[f^2(x)]dx) \prod_{j=1}^{n} K_j \bar{x}_j ^{2H_j - 1} \label{inequality}, 
	\end{equation}
 where
 \begin{align}
 K_j:=C_{H_j}^2 \Big\{4 (H_j-\tfrac{3}{2})^{-2}+2(H_j-\tfrac{1}{2})^{-1}(1-H_j)^{-1}\Big\}.\label{Lambda}
 \end{align}
	\end{theorem}
	
\begin{remark}
Note that we do not assume that  $f$ is adapted.
\end{remark}

\begin{remark}
Sufficient conditions (in terms of Malliavin derivatives) for the existence in $L^2(P)$ of the WIS integral with respect to 1-parameter fractional Brownian motion can be found in Bender \cite{Be}, Corollary 3.5.
For general $(\mathcal{S})^*$-valued integrands, conditions for the existence of the WIS-integral as an element of $(\mathcal{S})^*$ are given by Mishura \cite{M}, Theorem 2.3.1, using an inequality of Vaage.
\end{remark}

\dproof  By Theorem \ref{iso1} we have
\begin{align}
E[Y^2]&=E[(\int_{\R^n} f(x) dB_H(x))^2]= E[(\int_{\R^n} (Mf)(x)dB(x))^2]\nonumber\\
&= E[\int_{\R^n}  (Mf)^2(x) dx] = E[\int_{\R^n}f(x)M^2f(x)dx].
\end{align}
Also, by Lemma \ref{l1}

$$M^2f(x)= \int_{\mathbb{R}^n} \int_{\mathbb{R}^n} \prod_{j=1}^n C_{H_j}^2| y_j z_j| ^{H_j-\frac{3}{2}}  f(y+z+x)dy dz.$$

Hence, using the Cauchy-Schwartz inequality twice we obtain

\begin{align*}
E[Y^2]&=E[\int_{S} f(x) \int_{\mathbb{R}^n} \int_{\mathbb{R}^n} \prod_{j=1}^n C_{H_j}^2| y_j z_j| ^{H_j-\frac{3}{2}}  f(y+z+x)dy dz dx]\\
&=E[\int_{\R^n} \int_{\mathbb{R}^n} (\int_{S} \prod_{j=1}^n C_{H_j}^2| y_j z_j| ^{H_j-\frac{3}{2}}f(x)  f(y+z+x)dx)dy dz) ]\\
&=E[ \int_{\mathbb{R}^n} \int_{\mathbb{R}^n} \prod_{j=1}^n C_{H_j}^2| y_j z_j| ^{H_j-\frac{3}{2}} (\int_{S}f(x)  f(y+z+x)dx)dy dz]\\
&= \int_{\mathbb{R}^n} \int_{\mathbb{R}^n} \prod_{j=1}^n C_{H_j}^2| y_j z_j| ^{H_j-\frac{3}{2}} (\int_{S}E[f(x) f(y+z+x)]dx)dy dz\\
&\leq \int_{\mathbb{R}^n} \int_{\mathbb{R}^n} \prod_{j=1}^n C_{H_j}^2| y_j z_j| ^{H_j-\frac{3}{2}} (\int_{S}E[f^2(x)]^{\tfrac{1}{2}} E[f^2(y+z+x)]^{\tfrac{1}{2}}dx)dy dz\\
&\leq \int_{\mathbb{R}^n} \int_{\mathbb{R}^n} \prod_{j=1}^n C_{H_j}^2| y_j z_j| ^{H_j-\frac{3}{2}} (\int_{S}E[f^2(x)]dx \int_{S} E[f^2(y+z+x)]dx)^{\tfrac{1}{2}}dy dz\\
&=(\int_{S} E[f^2(x)]dx)^{\tfrac{1}{2}}  \int_{\mathbb{R}^n} \int_{\mathbb{R}^n} \prod_{j=1}^n C_{H_j}^2| y_j z_j| ^{H_j-\frac{3}{2}}  (\int_{S} E[f^2(y+z+x)]dx)^{\tfrac{1}{2}}dy dz\\
&\leq(\int_{S} E[f^2(x)]dx)^{\tfrac{1}{2}}  \int_{\{(y,z); y+z \in S\}} \prod_{j=1}^n C_{H_j}^2| y_j z_j| ^{H_j-\frac{3}{2}}  (\int_{S} E[f^2(y+z+x)]dx)^{\tfrac{1}{2}}dy dz\\
&\leq \int_{S} E[f^2(x)]dx  \int_{\{(y,z); y+z \in S\}} \prod_{j=1}^n C_{H_j}^2| y_j z_j| ^{H_j-\frac{3}{2}} dy dz\\
&\leq \int_{S} E[f^2(x)]dx  \prod_{j=1}^n C_{H_j}^2 \int_{0\leq y_j + z_j \leq \bar{x}_j} | y_j z_j| ^{H_j-\frac{3}{2}} dy_j dz_j.
\end{align*}
Here we have used that if $x+y+z \in S$ and $x \in S$, then $y+z \in S$, i.e. $0 \leq y_j + z_j \leq \bar{x}_j$ for all $j$.\\
It remains to estimate the integral
\begin{align}
J:= \int_{\{(y_j,z_j);0\leq y_j + z_j \leq \bar{x}_j\}} | y_j z_j| ^{H_j-\frac{3}{2}} dy_j dz_j \leq J_1 + J_2,
\end{align}
where
\begin{align}
    J_1= \int _{\{(y_j,z_j);|y_j|\leq \bar{x}_j,|z_j|\leq \bar{x}_j\}}| y_j z_j| ^{H_j-\frac{3}{2}} dy_j dz_j 
\end{align}
and
\begin{align}
    J_2= \int_{\{(y_j,z_j); |y_j|> \bar{x}_j,|z_j|> \bar{x}_j \text{ } \& \text{ } 0\leq y_j + z_j \leq \bar{x}_j\}}| y_j z_j| ^{H_j-\frac{3}{2}} dy_j dz_j 
\end{align}
By direct computation we get that
\begin{align}
    J_1 = 4 (H_j-\tfrac{3}{2})^{-2}\bar{x}_j^{2H_j -1}.
\end{align}

It is also easy to see by first integrating wrt. $dy_j$ and then using the mean value theorem, that
\begin{align}
    J_2 \leq 2 \bar{x}_j (H_j-\tfrac{1}{2})^{-1} \int_{\bar{x}_j}^{\infty}z_j^{2H_j - 3} dz_j \leq  2(H_j-\tfrac{1}{2})^{-1}(1-H_j)^{-1}\bar{x}_j^{2H_j-1}.
\end{align}

	We conclude that
	
	$$ E[Y^2]\leq \Big(\int_{S}E[f^2(x)]dx\Big) \prod_{j=1}^{n} C_{H_j}^2 \Big\{4 (H_j-\tfrac{3}{2})^{-2}+2(H_j-\tfrac{1}{2})^{-1}(1-H_j)^{-1}\Big\} \bar{x}_j ^{2H_j - 1},$$
	which completes the proof.
	
	\fproof
 \begin{remark}
 In this proof we have made no attempt to make the estimate as sharp as possible. It would be of interest to find the best constant in \eqref{inequality}.
 \end{remark}

\section{The fractional heat equation with additive noise}

In this section  we study the additive noise case, i.e. we study the equation \begin{align} \label{add 2}
 	 &\frac{\partial^{\alpha}}{\partial t^{\alpha}}Y(t,x)=\lambda \Delta Y(t,x)+\sigma  W_{H}(t,x);\; (t,x)\in (0,\infty)\times \mathbb{R}^{d}
   \end{align}
   with boundary conditions
\begin{align}\label{add01}
    Y(0,x)&=\delta_0(x)\text{ (the Dirac measure at  } x), \nonumber\\
    \lim_{x \rightarrow +/- \text{ }\infty}Y(t,x)&=0,
\end{align}
  \begin{definition}
   We say that a process $Y(t,x)=Y(t,x,\omega): [0,T] \times \mathbb{R}^d \times \Omega \mapsto \mathbb{R}$ is a solution of of the fractional stochastic heat equation \eqref{add 2} - \eqref{add01} if \\
(i) $\omega \mapsto Y(t,x,\omega) \in (\mathcal{S})^{*}$ for all $t,x$, \\
(ii) $(t,x) \mapsto Y(t,x,\omega) \in \mathcal{S}'$ for all $\omega \in \Omega$, and\\
(iii) the function $(t,x) \mapsto Y(t,x,\omega)$ satisfies the equation \eqref{add 2} - \eqref{add01} in the sense of distribution, for all $\omega$.
\end{definition}  
Our first main result in this section is the following:
 \begin{theorem} \label{add}
     The unique solution $Y(t,x) \in \mathcal{S}'$ of the fractional heat equation \eqref{add 2} - \eqref{add01} is given by
     \begin{align}\label{addsol}
     Y(t,x)&=(2\pi)^{-d} \int_{\mathbb{R}^d} e^{ixy}E_{\alpha}(-\lambda |y|^2 t^{\alpha})dy \nonumber\\
    &+ \sigma (2\pi)^{-d} \int_{0}^{t}(t-r)^{\alpha -1}\int_{0} ^x \Big(\int_{\mathbb{R}^{d}}e^{i(x-z)y}E_{\alpha}(-\lambda |y|^2 (t-r)^{\alpha})dy\Big) B_H(dr,dz).
    \end{align}
 \end{theorem}  
 \dproof
 In the following we let
\begin{align}
    L\{f(\cdot,x)\}(s,x)=\widetilde{f(\cdot,x)}(s,x)=\int_0^t  e^{-rs}f(r,x)dr;\quad s\geq 0
\end{align}
denote the Laplace transform of $f(t,x)$ as a function of $t$.\\

 In the following we let 
\begin{align}
    \L\{f(t,\cdot)\}(t,y)=\widehat{f(t,\cdot)}(t,y)
    \end{align}
    denote the (multi-dimensional) Laplace transform of $f(t,x)$ as a function of $x=(x_1, x_2, ... ,x_d).$ \\

First assume that $Y(t,x)$ is a solution of \eqref{add 2}.  We apply the Laplace transform $L$ 
to both sides of  \eqref{add 2} and obtain

\begin{equation}\label{step1a}
	s^{\alpha}\widetilde{Y}(s,x)-s^{\alpha-1}Y(0^{+},x)=\beta\widetilde{\Delta Y}(s,x)+\sigma \widetilde{W}_H(s,x).
	\end{equation}
Applying the multi-dimensional Laplace transform $\L$ to this equation  	 we get, since $\widehat{Y}(0,y)=1$, 
or, with $|y|^2=\sum_{j=1}^d y_j^2$,
   \begin{equation}\label{step2a}
   	\left( s^{\alpha}+\l|y|^2\right)\widehat{\widetilde{Y}}(s,y)=  s^{\alpha-1}\widehat{Y}(0^{+},y)+\sigma\widehat{\widetilde{W}_H}(s,y)
   \end{equation}
   Hence
   \begin{align}\label{step3}
       \widehat{\widetilde{Y}}(s,y)=\frac{s^{\alpha -1}}{s^{\alpha} + \l |y|^2}
       + \frac{\sigma \widehat{\widetilde{W}_H}(s,y)}{s^{\alpha}+\l |y|^2}.
   \end{align}
   Since the Laplace transform and the multi-dimensional Laplace   transform commute, this can be written
   \begin{align}\label{step4a}
       \widetilde{\widehat{Y}}(s,y)=\frac{s^{\alpha -1}}{s^{\alpha} + \l |y|^2}
       + \frac{\sigma \widetilde{\widehat{W}_H}(s,y)}{s^{\alpha}+\l |y|^2}.
   \end{align}
 Applying the inverse Laplace operator $L^{-1}$  to this equation we get
 \begin{align} 
       \widehat{Y}(t,y)&=L^{-1} \Big(\frac{s^{\alpha -1}}{s^{\alpha} + \l |y|^2}\Big)(t,y)
       + L^{-1}\Big(\frac{\sigma \widetilde{\widehat{W}_H}(s,y)}{s^{\alpha}+\l |y|^2}\Big)(t,y)\nonumber\\
       &=E_{\alpha,1}(-\l |y|^2 t^{\alpha}) + L^{-1}\Big(\frac{\sigma \widetilde{\widehat{W}_H}(s,y)}{s^{\alpha}+\l |y|^2}\Big)(t,y),\label{step 5}
   \end{align}
   where in general
   \begin{align} \label{add2.12}
   E_{\alpha,\beta}(z)=\sum_{k=0}^{\infty} \frac{z^{k}}{\Gamma(\alpha k + \beta)}    
   \end{align}
   is the Mittag-Leffler function.  \vskip 0.2cm
   It remains to find 
$L^{-1}\left(\frac{\sigma\widehat{\widetilde{W}_H}(s,y)}{s^{\alpha}+\l |y|^2}\right)$:\\
Recall that the convolution $f\ast g$ of two functions $f,g: [0,\infty) \mapsto \mathbb{R}$  is defined by
\begin{align}
(f \ast g)(t)=\int_0^t f(t-r)g(r) dr; \quad t \geq 0.
\end{align}
The convolution rule for Laplace transform states that $$L\left( \int_{0}^{t}f(t-r)g(r)dr\right) (s)=Lf(s)Lg(s),$$ 
or 
\begin{equation}\label{addititive 12}
	\int_{0}^{t}f(t-w)g(w)dw=L^{-1}\left( Lf(s)Lg(s)\right) (t).\\
\end{equation}
Now (see \cite{I},p.312)
\begin{align}
	L^{-1}\left(\frac{1}{s^{\alpha}+\l |y|^2} \right) (t)&=t^{\alpha-1}E_{\alpha,\alpha}(-\l t^{\alpha}|y|^2)\nonumber\\
	&=\sum_{k=0}^{\infty}\frac{t^{\alpha-1}(-\l t^{\alpha}|y|^2)^{k}}{\Gamma(\alpha k+\alpha)}\nonumber\\
	&=\sum_{k=0}^{\infty}\frac{(-\l |y|^2)^{k}t^{\alpha(k+1)-1}}{\Gamma(\alpha(k+1))}\nonumber\\
	&=\sum_{k=0}^{\infty}\frac{(-\l t^{\alpha}|y|^2)^{k}t^{\alpha -1}}{\Gamma(\alpha (k+1))}	\nonumber\\
	&=: \Lambda(t,y).\label{addL}
\end{align}
In other words,
\begin{equation}
	\frac{\sigma}{s^{\alpha}+\l |y|^2}=\sigma L \Lambda(s,y),
\end{equation}
Combining with \eqref{addititive 12} we get 
\begin{align}
	L^{-1}\left( \frac{\sigma}{s^{\alpha}+\l |y|^2} \widehat{\widetilde{W}} (s,y)\right) (t)&=L^{-1}\left( L\left( \sigma\Lambda(s,y)\right) \widetilde{\widehat{W}}(s,y)\right) (t)\\
	&=\sigma \int_{0}^{t}\Lambda(t-r,y)\widehat{W}_H(r,y)dr.
\end{align}
Substituting this into \eqref{step 5} we get 
\begin{equation}
	\widehat{Y}(t,y)=E_{\alpha,1}\left( -\l t^{\alpha}|y|^2\right)+\sigma\int_{0}^{t}\Lambda(t-r,y)\widehat{W}_H(r,y)dr.
\end{equation}
Taking inverse multidimensional Laplace  transform we end up with
\begin{equation}\label{add2.18a}
 Y(t,x)=\L^{-1}\left( E_{\alpha,1}\left(-\l t^{\alpha}|y|^2 
 \right)\right)(x)+\sigma \L^{-1}\left(\int_{0}^{t}\Lambda(t-r,y)\widehat{W}_H(r,y)dr\right)(x).
\end{equation}
 
Now, we use that 
\begin{align}
\L\{\int_{0}^{x}f(x-z)g(z)dz\}(y)= \L\{f\}(y)\L\{g\}(y),\quad x \in \mathbb{R}_{+}^d
\end{align}
or
\begin{equation}
	\int_{0}^{x}f(x-z)g(z)dz=\L^{-1}\Big( \L\{f\}(y)\L\{g\}(y)\Big)(x)\quad x \in \mathbb{R}_{+}^d,
\end{equation}
 where we are using the notation
 $$\int_0^x \varphi(z)dz=\int_0^{x_d} \int_0^{x_{d-1}} ... \int_0^{x_1} \varphi(z_1,z_2, ... , z_d) dz_1 dz_2 ... dz_1 $$ 
 when $x=(x_1, x_2, ... , x_d), z=(z_1, z_2, ... ,z_d)$,

This gives,

\begin{align*}
	\L^{-1}\left( \int_{0}^{t}\Lambda(t-r,y)\widehat{W}(r,y)dr\right) (x)&=\int_{0}^{t}\L^{-1}\left( \Lambda(t-r,y)\widehat{W}_H(r,y)\right) (x)dr\\
	&=\int_{0}^{t}\L^{-1}\Big( \L\left( \L^{-1}\Lambda(t-r,y)\right) (y)\L W_H(r,x)(y)\Big) (x)dr\\
	&=\int^x _{0}\int_{0}^{x}\left( \L^{-1}\Lambda(t-r,y)(x-z)\right) W_H(r,z)dzdr\\
	&=\int_{0}^{t}\int_{0}^{x}\left( (2\pi)^{-d}\int_{\mathbb{R}^{d}}e^{i(x-z)y}\Lambda(t-r,y)dy\right) W_H(r,z)dzdr\\
	&=(2\pi)^{-d}\int_{0}^{t}\int_{0}^x \left( \int_{\mathbb{R}^{d}}e^{i(x-z)y}\Lambda(t-r,y)dy\right) B_H(dr,dz).
\end{align*}
Here we have used that the function
\begin{align}
   y \mapsto  \Lambda(t,y)=\sum_{k=0}^{\infty}\frac{(-\l t^{\alpha}y^2)^{k}t^{\alpha -1}}{\Gamma(\alpha (k+1))}; \quad y \in \mathbb{R}^d
\end{align}
can be extended to an entire complex function in the complex space $\mathbb{C}^d$, and hence its inverse Laplace transform coincides with the inverse Fourer transform.\\
Combining this with \eqref{add2.18}, \eqref{add2.12} and \eqref{addLambda} we get
\begin{align}\label{addexp1a}
    Y(t,x)&= \L^{-1}(\sum_{k=0}^{\infty} \frac{(- \l t^{\alpha} |y|^2)^k}{\Gamma(\alpha k +1)})+\sigma (2\pi)^{-d} \int_{0}^{t}\int_{0}^x \left(\int_{\mathbb{R}^{d}}e^{i(x-z)y}\Lambda(t-r,y)dy\right) B_H(dr,dz)\nonumber\\
    &=(2\pi)^{-d} \int_{\mathbb{R}^d} e^{ixy}\sum_{k=0}^{\infty} \frac{(- \l t^{\alpha} |y|^2)^k}{\Gamma(\alpha k +1)}dy \nonumber\\
    &+ \sigma (2\pi)^{-d} \int_{0}^{t}(t-r)^{\alpha -1}\int_{0} ^x \left(\int_{\mathbb{R}^{d}}e^{i(x-z)y}\sum_{k=0}^{\infty}\frac{(-\l (t-r)^{\alpha}|y|^2)^{k}}{\Gamma(\alpha (k+1))}dy\right) B_H(dr,dz) \\
    &=(2\pi)^{-d} \int_{\mathbb{R}^d} e^{ixy}E_{\alpha}(-\lambda |y|^2 t^{\alpha})dy \nonumber\\
    &+ \sigma (2\pi)^{-d} \int_{0}^{t}(t-r)^{\alpha -1}\int_{0} ^x \Big(\int_{\mathbb{R}^{d}}e^{i(x-z)y}E_{\alpha}(-\lambda |y|^2 (t-r)^{\alpha})dy\Big) B_H(dr,dz)
\end{align}
This proves uniqueness and also that the unique solution (if it exists) is given by \eqref{addsol}.\\

b) Next, define $Y(t,x)$ by \eqref{addexp1}. Then we see that we can prove that $Y(t,x)$ satisfies \eqref{add} by reversing the argument above.
We skip the details.
\fproof

\fproof

\subsection{The classical case ($\alpha$ = 1)}
It is interesting to compare the above result with the classical case when $\alpha$=1:\\
If $\alpha=1$, we get $Y(t,x)=I_{0}(x)+I_{1}(t,x)$, where
\begin{align*}
	&I_{0}(x)=(2\pi)^{-d}\int_{\mathbb{R}^{d}}e^{ixy}\sum_{k=0}^{\infty}\frac{\left( -\l t |y|^{2}\right) ^{k}}{k!}dy=(2\pi)^{-d} \int_{\mathbb{R}^d} e^{ixy}E_{\alpha}(-\lambda |y|^2 t^{\alpha})dy, \\
&I_1(t,x)= \sigma (2\pi)^{-d} \int_{0}^{t}\int_{0}^x \Big(\int_{\mathbb{R}^{d}}e^{i(x-z)y}E_{\alpha}(-\lambda |y|^2 (t-r)^{\alpha})\Big) B_H(dr,dz)
\end{align*}
Using that that $\Gamma(k+1)=k!$ 
and the Taylor expansion of the exponential function, we get
\begin{align*}
 I_{0}&=(2\pi)^{-d}\int_{\mathbb{R}^{d}}e^{ixy}e^{-\l t |y|^{2}} dy\\
 &=(2\pi)^{-d}\left( \frac{\pi}{\l t}\right) ^{\frac{d}{2}}e^{-\frac{|x|^{2}}{4\l t}}\\
 &=(4 \pi \l t) ^{-\frac{d}{2}} e^{-\frac{|x|^{2}}{4\l t}},	
\end{align*}
where we used the general formula
\begin{equation}
\int_{\mathbb{R}^{d}}e^{-\left( a |y|^{2}+2by\right) }dy=\left(\frac{\pi}{a} \right)^{\frac{d}{2}} e^{\frac{b^{2}}{a}};\;a>0;\;b\in\mathbb{C}^d.\label{exp}
\end{equation}\\
Similarly,
\begin{align*}
	I_{1}(t,x)&=\sigma(2\pi)^{-d}\int_{0}^{t}\int_{0} ^x \int_{\mathbb{R}^{d}}e^{i(x-z)y}\sum_{k=0}^{\infty}\frac{\left(-\l (t-r)|y|^{2} \right)^{k} }{k!}dyB_H(dr,dz)\\
	&=\sigma(2\pi)^{-d}\int_{0}^{t}\int_{0}^x \left( \frac{\pi}{\l (t-r)}\right) ^{\frac{d}{2}}e^{-\frac{|x-z|^{2}}{4\l (t-r) }}B_H(dr,dz)\\
	&=\sigma (4 \pi \l)^{-\tfrac{d}{2}}\int_0^t \int_{0}^x  (t-r)^{-\tfrac{d}{2}} e^{-\frac{|x-z|^2}{4 \l (t-r)}} B_H(dr,dz).
\end{align*}
Summarizing the above, we get, for $\alpha=1$, 
\begin{align}
	Y(t,x)&=(4\pi \l t)^{-\frac{d}{2}}e^{-\frac{|x|^{2}}{4\l t}}\nonumber\\
	&+\sigma(4\pi\l)^{-\frac{d}{2}}\int_{0}^{t}\int_{0}^x (t-r)^{-\frac{d}{2}}e^{-\frac{|x-z|^{2}}{ 4\l (t-r)}}B_H(dr,dz)\;\; \label{alpha1}
\end{align}
This is in agreement with a well-known classical result. See e.g. Section 4.1 in \cite{Hu}.

\subsection{When is $Y(t,x)$ a mild solution?}

It was pointed out already in 1984 by John Walsh \cite{W} that (classical) SPDEs driven by time-space white noise $W(t,x); (t,x) \in [0,\infty) \times \R^d$ may have only distribution valued solutions if $d \geq 2$. 
Indeed, the solution $Y(t,x)$ that we found in the previous section is in general distribution valued. But in some cases the solution can be represented as an element of $L^2(\P)$. Following Y. Hu \cite{Hu} we make the following definition:

\begin{definition}
The solution $Y(t,x)$ is called \emph{mild} if $Y(t,x) \in L^2(\P)$ \\for all $t>0, x >0$.
\end{definition}

We now ask the following:
\begin{problem}
   For what values of $\alpha \in (0,2)$ and what dimensions $d=1,2, ...$ is $Y(t,x)$ mild?
      \end{problem}

By the Ito  isometry  we have 
\begin{equation}
    \E\Big[ Y^{2}(t,x)\Big]=J_{1}+J_{2},
\end{equation}
where
\small
 \begin{align}
 J_1&=(2\pi)^{-2d} \Big(\int_{\mathbb{R}^d} e^{ixy}\sum_{k=0}^{\infty} \frac{(- \l t^{\alpha} |y|^2)^k}{\Gamma(\alpha k +1)}dy\Big)^2 \nonumber\\
   J_2 &=\E\Big[ \Big(\sigma (2\pi)^{-d} \int_{0}^{t}\int_{0} ^x(t-r)^{\alpha -1} \left(\int_{\mathbb{R}^{d}}e^{i(x-z)y}\sum_{k=0}^{\infty}\frac{(-\l (t-r)^{\alpha}|y|^2)^{k}}{\Gamma(\alpha (k+1))}dy\right) B_H(dr,dz) \Big)^2\Big]\nonumber\\
   &\sim   \int_{0}^{t}\int_{0} ^x(t-r)^{2\alpha -2} \left(\int_{\mathbb{R}^{d}}e^{i(x-z)y}\sum_{k=0}^{\infty}\frac{(-\l (t-r)^{\alpha}|y|^2)^{k}}{\Gamma(\alpha (k+1))}dy\right)^2 dr dz ,
 \end{align}
 where the notation
 $f\sim g$ means that there exist functions $c(t,x)>0, C(t,x)>0$  such that $c(t,x) f(t,x) \leq g(t,x) \leq C(t,x) f(t,x).$

This is the same estimate as in Theorem in \cite{MO}. Therefore the proof of that theorem applies here, and we can conclude the following:
 
  \begin{theorem} \label{mild1}
 Let $Y(t,x)$ be the solution of the $\alpha$-fractional stochastic heat equation with additive noise.
Then the following holds for all $H=(H_1,H_2, ... ,H_d) \in (\tfrac{1}{2},1)^d$:
\begin{itemize}
\item
a) If $\alpha = 1$, then $Y(t,x)$ is mild if and only if $d=1$.
\item
b) If $\alpha > 1$ then $Y(t,x)$ is mild if $d=1$ or $d=2$.
\item
c) If $\alpha < 1$ then $Y(t,x)$ is not mild for any  $d.$
\end{itemize}
 \end{theorem} 

\begin{remark}
\begin{itemize}
\item
See Y. Hu \cite{Hu}, Proposition 4.1 for a generalisation of the above result in the case $\alpha=1$.
\item
In the cases $\alpha > 1, d \geq 3$ we do not know if the solution $Y(t,x)$ is mild or not. This is a topic for future research.
\end{itemize}
\end{remark}

\section{The fractional heat equation with multiplicative noise}
In this section we study the multiplicative noise case, i.e. we study the equation
\begin{align}\label{mult}
&\frac{\partial^{\alpha}}{\partial t^{\alpha}}Y(t,x)=\lambda \Delta Y(t,x)+\sigma Y(t,x) \diamond W_{H}(t,x);\; (t,x)\in (0,\infty)\times \mathbb{R}^{d}
\end{align}
   with boundary conditions
\begin{align}\label{mult0}
    Y(0,x)&=\delta_0(x)\text{ (the Dirac measure at  } x), \nonumber\\
    \lim_{x \rightarrow +/- \text{ }\infty}Y(t,x)&=0,
\end{align}
 We first clarify what we mean by a solution:
\begin{definition}
We say that a process $Y(t,x)=Y(t,x,\omega): [0,T] \times \mathbb{R}^d \times \Omega \mapsto \mathbb{R}$ is a solution of of the fractional stochastic heat equation \eqref{mult} - \eqref{mult0} if \\
(i) $\omega \mapsto Y(t,x,\omega) \in (\mathcal{S})^{*}$ for all $t,x$, \\
(ii) $(t,x) \mapsto Y(t,x,\omega) \in \mathcal{S}'$ for all $\omega \in \Omega$, and\\
(iii) the function $(t,x) \mapsto Y(t,x,\omega)$ satisfies the equation \eqref{add 2} - \eqref{1.5} in the sense of distribution, for all $\omega$.
\end{definition}

Our first main result in this section is the following:
\begin{theorem} \label{mult1}
A process $Y(t,x) $ is a solution of the of the fractional stochastic heat equation \eqref{mult} - \eqref{mult0} if and only if $Y(t,x)$ is a solution of the following time-space stochastic Volterra equation:
\begin{align} \label{vol}
Y(t,x)&=f(t,x) + \int_0^t \int_{0}^x  g(t-r,x-z) Y(r,z) B_H(dr,dz),
\end{align}
where
\small{
\begin{align}
f(t,x)=(2\pi)^{-d} \int_{\mathbb{R}^d} e^{ixy} E_{\alpha}(- \l t^{\alpha} |y|^2) dy
=(2\pi)^{-d} \int_{\mathbb{R}^d} e^{ixy}\sum_{k=0}^{\infty} \frac{(- \l t^{\alpha} |y|^2)^k}{\Gamma(\alpha k +1)}dy, 
\end{align}
}
and
\small{
\begin{align}
    &g(t-r,x-z)= \sigma (2\pi)^{-d} (t-r)^{\alpha -1}\int_{\mathbb{R}^{d}}e^{i(x-z)y} E_{\alpha,\alpha}(-\l (t-r)^{\alpha}|y|^2) dy
\end{align}
}
where $|y|^2=y^2=\sum_{j=1}^d y_j^2.$ 
\end{theorem}
\begin{remark}
We say that $Y(t,x)$ solves the fractional Volterra equation \eqref{vol} if $Y(t,x)$ is an $(\mathcal{S})^{*}$-valued solution of the equation
\begin{align} \label{3.1a}
Y(t,x)&=f(t,x) + \int_0^t \int_{0} ^x g(t-r,x-z) Y(r,z) \diamond W_{H}(r,z) dr dz,
\end{align}
\end{remark}

Proof of Theorem \ref{mult1}:\\
a) We proceed as in the proof of Theorem \ref{add}:\\
First assume that $Y(t,x)$ is a solution of \eqref{mult0}.  We apply the Laplace transform $L$ 
to both sides of  \eqref{add 2} and obtain

\begin{equation}\label{step1}
	s^{\alpha}\widetilde{Y}(s,x)-s^{\alpha-1}Y(0^{+},x)=\beta\widetilde{\Delta Y}(s,x)+\sigma \widetilde{Y W}_H(s,x).
	\end{equation}
Applying the multi-dimensional Laplace transform $\L$ to this equation  	 we get, since $\widehat{Y}(0,y)=1$, 
or, with $|y|^2=\sum_{j=1}^d y_j^2$,
   \begin{equation}\label{step2}
   	\left( s^{\alpha}+\l|y|^2\right)\widehat{\widetilde{Y}}(s,y)=  s^{\alpha-1}\widehat{Y}(0^{+},y)+\sigma\widehat{\widetilde{Y W}_H}(s,y)
   \end{equation}
   Hence
   \begin{align}\label{step3a}
       \widehat{\widetilde{Y}}(s,y)=\frac{s^{\alpha -1}}{s^{\alpha} + \l |y|^2}
       + \frac{\sigma \widehat{\widetilde{Y W}_H}(s,y)}{s^{\alpha}+\l |y|^2}.
   \end{align}
   Since the Laplace transform and the multi-dimensional Laplace   transform commute, this can be written
   \begin{align}\label{step4}
       \widetilde{\widehat{Y}}(s,y)=\frac{s^{\alpha -1}}{s^{\alpha} + \l |y|^2}
       + \frac{\sigma \widetilde{\widehat{Y W}_H}(s,y)}{s^{\alpha}+\l |y|^2}.
   \end{align}
 Applying the inverse Laplace operator $L^{-1}$  to this equation we get
 \begin{align} 
       \widehat{Y}(t,y)&=L^{-1} \Big(\frac{s^{\alpha -1}}{s^{\alpha} + \l |y|^2}\Big)(t,y)
       + L^{-1}\Big(\frac{\sigma \widetilde{\widehat{Y W}_H}(s,y)}{s^{\alpha}+\l |y|^2}\Big)(t,y)\nonumber\\
       &=E_{\alpha,1}(-\l |y|^2 t^{\alpha}) + L^{-1}\Big(\frac{\sigma \widetilde{\widehat{Y W}_H}(s,y)}{s^{\alpha}+\l |y|^2}\Big)(t,y).\label{step5a}
   \end{align}
    \vskip 0.2cm
   It remains to find 
$L^{-1}\left(\frac{\sigma\widehat{\widetilde{Y W}_H}(s,y)}{s^{\alpha}+\l |y|^2}\right)$:\\
Recall that the convolution $f\ast g$ of two functions $f,g: [0,\infty) \mapsto \mathbb{R}$  is defined by
\begin{align}
(f \ast g)(t)=\int_0^t f(t-r)g(r) dr; \quad t \geq 0.
\end{align}
The convolution rule for Laplace transform states that $$L\left( \int_{0}^{t}f(t-r)g(r)dr\right) (s)=Lf(s)Lg(s),$$ 
or 
\begin{equation}\label{add12}
	\int_{0}^{t}f(t-w)g(w)dw=L^{-1}\left( Lf(s)Lg(s)\right) (t).\\
\end{equation}
Now (see \cite{I},p.312)
\begin{align}
	L^{-1}\left(\frac{1}{s^{\alpha}+\l |y|^2} \right) (t)&=t^{\alpha-1}E_{\alpha,\alpha}(-\l t^{\alpha}|y|^2)\nonumber\\
	&=\sum_{k=0}^{\infty}\frac{t^{\alpha-1}(-\l t^{\alpha}|y|^2)^{k}}{\Gamma(\alpha k+\alpha)}\nonumber\\
	&=\sum_{k=0}^{\infty}\frac{(-\l |y|^2)^{k}t^{\alpha(k+1)-1}}{\Gamma(\alpha(k+1))}\nonumber\\
	&=\sum_{k=0}^{\infty}\frac{(-\l t^{\alpha}|y|^2)^{k}t^{\alpha -1}}{\Gamma(\alpha (k+1))}	\nonumber\\
	&=: \Lambda(t,y).\label{addLambda}
\end{align}
In other words,
\begin{equation}
	\frac{\sigma}{s^{\alpha}+\l |y|^2}=\sigma L \Lambda(s,y),
\end{equation}
Combining with \eqref{addititive 12} we get 
\begin{align}
	L^{-1}\left( \frac{\sigma}{s^{\alpha}+\l |y|^2} \widehat{\widetilde{Y W}} (s,y)\right) (t)&=L^{-1}\left( L\left( \sigma\Lambda(s,y)\right) \widetilde{\widehat{Y W}}(s,y)\right) (t)\\
	&=\sigma \int_{0}^{t}\Lambda(t-r,y)\widehat{Y W}_H(r,y)dr.
\end{align}
Substituting this into \eqref{step 5} we get 
\begin{equation}
	\widehat{Y}(t,y)=E_{\alpha,1}\left( -\l t^{\alpha}|y|^2\right)+\sigma\int_{0}^{t}\Lambda(t-r,y)\widehat{Y W}_H(r,y)dr.
\end{equation}
Taking inverse multidimensional Laplace  transform we end up with
\begin{equation}\label{add2.18}
 Y(t,x)=\L^{-1}\left( E_{\alpha,1}\left(-\l t^{\alpha}|y|^2 
 \right)\right)(x)+\sigma \L^{-1}\left(\int_{0}^{t}\Lambda(t-r,y)\widehat{Y W}_H(r,y)dr\right)(x).
\end{equation}
 
Now, we use that 
\begin{align}
\L\{\int_{0}^{x}f(x-z)g(z)dz\}(y)= \L\{f\}(y)\L\{g\}(y),\quad x \in \mathbb{R}_{+}^d
\end{align}
or
\begin{equation}
	\int_{0}^{x}f(x-z)g(z)dz=\L^{-1}\Big( \L\{f\}(y)\L\{g\}(y)\Big)(x)\quad x \in \mathbb{R}_{+}^d,
\end{equation}
 where we are using the notation
 $$\int_0^x \varphi(z)dz=\int_0^{x_d} \int_0^{x_{d-1}} ... \int_0^{x_1} \varphi(z_1,z_2, ... , z_d) dz_1 dz_2 ... dz_1 $$ 
 when $x=(x_1, x_2, ... , x_d), z=(z_1, z_2, ... ,z_d)$,

This gives,

\begin{align*}
	\L^{-1}\left( \int_{0}^{t}\Lambda(t-r,y)\widehat{Y W}(r,y)dr\right) (x)&=\int_{0}^{t}\L^{-1}\left( \Lambda(t-r,y)\widehat{Y W}_H(r,y)\right) (x)dr\\
	&=\int_{0}^{t}\L^{-1}\Big( \L\left( \L^{-1}\Lambda(t-r,y)\right) (y)\L \{Y W_H\}(r,x)(y)\Big) (x)dr\\
	&=\int^x _{0}\int_{0}^{x}\left( \L^{-1}\Lambda(t-r,y)(x-z)\right) Y(r,z)W_H(r,z)dzdr\\
	&=\int_{0}^{t}\int_{0}^{x}\left( (2\pi)^{-d}\int_{\mathbb{R}^{d}}e^{i(x-z)y}\Lambda(t-r,y)dy\right) Y(r,z) W_H(r,z)dzdr\\
	&=(2\pi)^{-d}\int_{0}^{t}\int_{0}^x \left( \int_{\mathbb{R}^{d}}e^{i(x-z)y}\Lambda(t-r,y)dy\right) Y(r,z) B_H(dr,dz).
\end{align*}
Here we have used that the function
\begin{align}
   y \mapsto  \Lambda(t,y)=\sum_{k=0}^{\infty}\frac{(-\l t^{\alpha}y^2)^{k}t^{\alpha -1}}{\Gamma(\alpha (k+1))}
\end{align}
can be extended to an entire complex function in the complex plane, and hence its inverse Laplace transform coincides with the inverse Fourer transform.\\
Combining this with \eqref{add2.18}, \eqref{add2.12} and \eqref{addLambda} we get
\begin{align}\label{addexp1}
    Y(t,x)&= \L^{-1}(\sum_{k=0}^{\infty} \frac{(- \l t^{\alpha} |y|^2)^k}{\Gamma(\alpha k +1)})+\sigma (2\pi)^{-d} \int_{0}^{t}\int_{0}^x \left(\int_{\mathbb{R}^{d}}e^{i(x-z)y}\Lambda(t-r,y)dy\right) Y(r,z) B_H(dr,dz)\nonumber\\
    &=(2\pi)^{-d} \int_{\mathbb{R}^d} e^{ixy}\sum_{k=0}^{\infty} \frac{(- \l t^{\alpha} |y|^2)^k}{\Gamma(\alpha k +1)}dy \nonumber\\
    &+ \sigma (2\pi)^{-d} \int_{0}^{t}(t-r)^{\alpha -1}\int_{0} ^x \left(\int_{\mathbb{R}^{d}}e^{i(x-z)y}\sum_{k=0}^{\infty}\frac{(-\l (t-r)^{\alpha}|y|^2)^{k}}{\Gamma(\alpha (k+1))}dy\right) Y(r,z) B_H(dr,dz) \\
    &=(2\pi)^{-d} \int_{\mathbb{R}^d} e^{ixy}E_{\alpha}(-\lambda |y|^2 t^{\alpha})dy \nonumber\\
    &+ \sigma (2\pi)^{-d} \int_{0}^{t}(t-r)^{\alpha -1}\int_{0} ^x \Big(\int_{\mathbb{R}^{d}}e^{i(x-z)y}E_{\alpha}(-\lambda |y|^2 (t-r)^{\alpha})dy\Big) Y(r,z) B_H(dr,dz)
\end{align}
This proves uniqueness and also that the unique solution (if it exists) is given by \eqref{addsol}.\\

b) Next, define $Y(t,x)$ by \eqref{addexp1}. Then we see that we can prove that $Y(t,x)$ satisfies \eqref{add} by reversing the argument above.
We skip the details.
\fproof

\section{A solution formula of the Volterra equation}

In the is section we find a formula for the Volterra eqiation \eqref{vol} and we use it to study the properties of the solution.\\

We will need the following general result about stochastic linear Volterra equations, based on \cite{OZ}:
\begin{theorem}\label{th8.1}
Let $F(u):  [0,\infty)^n \mapsto \mathbb{R}$ and $G(u,w): [0,\infty)^n \times [0,\infty)^n \mapsto \mathbb{R}$ be locally bounded deterministic functions such that
$F \in L^2(\mathbb{R}_{+}^n)$ and 
\begin{align}
\| G(u,\cdot)\|^2 &:= \| G(u,\cdot)\|_{L^{2}([0,u])} ^2:=   \int_0^u G^2(u,w) dw \nonumber\\
&:=\int_0^{u_1} \int_0^{u_2} ... \int_0^{u_n} G^2(u_1, u_2, ... ,u_n,w_1,w_2 ... w_n) dw_n ... dw_2 dw_1 < \infty
\end{align}
for all $u=(u_1, u_2, ... , u_n) \in \mathbb{R}_{+}^n$.
Then there is a unique solution $Y(u) \in L^2(P)$ of the stochastic Volterra equation
\begin{align}
    Y(u)= F(u) + \int_0^u G(u,w) Y(w) B_H(dw).
\end{align}
The solution is given by
\begin{align}
    Y(u)= F(u) + \int_0^u H(u,w) F(w) B_H(dw).
\end{align}
where
\begin{align}
    H(u,w)= \sum_{m=1}^{\infty} G_m(u,w)
\end{align}
with $G_m$ given inductively by
\begin{align}
G_{m+1}(u,w)&=\int_w ^u G_m(u,v)G(v,w) B_H(dv);\quad m=1,2, ...\nonumber\\
G_1(u,w)&=G(u,w).
\end{align}

\end{theorem}

\dproof
The proof is based on a suitably adapted modification of the proof of Theorem 3.1 in \cite{OZ}. We give an outline:\\

By induction we see that
\begin{align}
    G_m(u,w)=\underset{S(u,w)}{\int ... \int } \prod_{k=0}^{m-1} G(z_k,z_{k+1}) B_H(dz_1) B_H(dz_2) ... B_H(dz_{m-1}),
\end{align}
where $z_0=u, z_n=w$ and the integration is taken over the set
\begin{align}
S(u,w):= \Big\{(z_1,z_2, ... ,z_{m-1}) \in (\mathbb{R}_{+}^n)^{m-1}; w\leq z_{m-1}\leq z_{m-2}\leq  ...  \leq z_{1} \leq u \Big\}
\end{align}
Here we use the notation $$ x=(x_1, x_2, ... ,x_n) \leq y= (y_1,y_2, ... ,y_n) \text{ iff } x_j \leq y_j; \text {for all } j=1,2, ... ,n.$$
By applying Theorem \ref{th8.1} repeatedly we get that
\begin{align}
E[G_m(u,w)^2] &\leq A^{m-1}  \int_{S(u,w)} \prod_{k=0}^{m-1} G^2(z_k,k_{k+1})dz_1 dz_2 ... dz_{m-1}\nonumber\\
&\leq A^{m-1} \|G(u,\cdot)\| ^{2m} \frac{Vol(S(u,w))}{((m-1)!)^n}
\end{align} 
where
\begin{align}
A=\prod_{j=1}^n K_j \bar{x}_j^{2H_j - 1} \text{ with } K_j:=C_{H_j}^2 \Big\{4 (H_j-\tfrac{3}{2})^{-2}+2(H_j-\tfrac{1}{2})^{-1}(1-H_j)^{-1}\Big\}.\label{Lambda1}
 \end{align}
 Hence

$H(u,w)= \sum_{m=0}^{\infty} G_m(u,w)$ converges in $L^2(P)$ for all $u,w$.\\
By Theorem \ref{th8.1} it follows that
\begin{align}
 Y(u)=F(u) + \int_0 ^{u} H(u,w) F(w) B_H(dw)   \in L^2 (P) \text{ for all }u.
\end{align}
\fproof

We now apply this general result to our situation and obtain the following:

\begin{theorem}{(Solution of the fractional Volterra equation)}\label{vol3}
\newline
As before put
\small{
\begin{align}
f(t,x)=(2\pi)^{-d} \int_{\mathbb{R}^d} e^{ixy} E_{\alpha}(- \l t^{\alpha} |y|^2) dy
=(2\pi)^{-d} \int_{\mathbb{R}^d} e^{ixy}\sum_{k=0}^{\infty} \frac{(- \l t^{\alpha} |y|^2)^k}{\Gamma(\alpha k +1)}dy, 
\end{align}
}
and

\begin{align}
    &g(t-r,x-z)= \sigma (2\pi)^{-d} (t-r)^{\alpha -1}\int_{\mathbb{R}^{d}}e^{i(x-z)y} E_{\alpha,\alpha}(-\l (t-r)^{\alpha}|y|^2) dy  
\end{align}
 The solution of the fractional Volterra equation (\eqref{vol}): 
 \begin{align} 
Y(t,x)&=f(t,x) + \int_0^t \int_{0}^x  g(t-s,x-z) Y(s,z) B_H(ds,dz),
\end{align} 
 is
 \begin{align}
     Y(t,x)= f(t,x) + \int_0^t \int_0^x h(t,x,s,z) f(s,z) B_H(ds,dz),
 \end{align}
 where 
 \begin{align}
     h(t,x,s,z)=\sum_{m=1}^{\infty} g_m(t,x,s,z).
 \end{align}
with $g_m$ given inductively by
\begin{align}
g_{m}(t,x,s,z)&=\int_0 ^t\int_0^x \Big(\int_0^s \int_0^z g_{m-1}(t,x,r,w)g_1(s,z,r,w) B_H(dr,dw)\Big)B_H(ds,dz);\quad m=2,3, ...\nonumber\\
g_1(t,x,s,z)&=g(t-s,x-z),
\end{align}
i.e.
\begin{align*}
    &g_m(t,x,s,z)\nonumber\\
    &=\underset{S((t,x),(s,z))}{\int ... \int } 
    \prod_{k=0}^{m-1} g(s_k,z_k,s_{k+1},z_{k+1}) B_H(ds_1,dz_1) B_H(ds_2,dz_2) ...  B_H(ds_{m-1},dz_{m-1}).
\end{align*}
Here $(s_0,z_0)=(t,x), (s_{m},z_{m})=(s,z)$ and the integration is taken over the set
\begin{align}
S((t,x),(s,z))&:= \Big\{(s_1,z_1,s_2,z_2, ... , s_{m-1},z_{m-1})\in (\mathbb{R}_{+}^{1+d})^{m-1};\nonumber\\
&s\leq s_{m-1} \leq s_{m-2} \leq ... \leq t\quad  \& \quad z \leq z_{m-1}\leq ... \leq z_1 \leq x \quad  \Big\}.
\end{align}
\end{theorem}

\subsection{When is the solution mild?}
We now apply the result above to the Volterra equation \eqref{vol}:\\
In this case we have $u=(t,x),w=(r,z)$  and
\begin{align}
    G(u,w)= \sigma (2\pi)^{-d} (t - r)^{\alpha -1} \int_{\mathbb{R}^d} e^{i(x-z)y} E_{\alpha,\alpha}( -\l (t-r)^{\alpha}|y|^2)dy. 
\end{align}
Hence
\begin{align}
    \|G(u,\cdot)\|^2 &= \sigma^2 (2\pi)^{-2d} \int_0^t \int_0^x (t-r)^{2\alpha -2}\Big( \int_{\mathbb{R}^d} e^{i(x-z)y} E_{\alpha,\alpha}( -\l (t-r)^{\alpha}|y|^2)dy\Big)^2 dr dz\nonumber\\
& \leq \sigma^2 (2\pi)^{-2d} \int_0^t \int_0^x (t-r)^{2\alpha -2}\Big( \int_{\mathbb{R}^d} | E_{\alpha,\alpha}( -\l (t-r)^{\alpha}|y|^2)|dy\Big)^2 dr dz.\nonumber\\ \label{11.38}
\end{align}
We split into two cases:\\
1) Assume $\alpha \geq 1$.\\
Then by Jia et al \cite{J} we have
\begin{align}
    |E_{\alpha,\alpha}(-\l (t-r)^{\alpha} |y|^2)| \leq e^{-\l (t-r)^{\alpha} |y|^2}. 
\end{align}
Hence
\begin{align}
    \|G(u,\cdot) \| ^2  \leq \sigma^2 (2\pi)^{-2d} \int_0^t \int_0^x (t-r)^{2\alpha -2}\Big( \int_{\mathbb{R}^d}  e^{-\l (t-r)^{\alpha} |y|^2}dy\Big)^2 dr dz.
\end{align}
If we introduce the change of variable $v=(\l (t-r)^{\alpha})^{\tfrac{1}{2}} y$ we get
\begin{align}
\|G(u,\cdot)\|^2 &\leq \int_0^t (t-r)^{2\alpha -2} ( \l (t-r)^{\alpha})^{-d} (\int _{\mathbb{R}^d  }  e^{-v^2} dv)^2 \prod_{j=1}^d x_j dr\nonumber\\ 
 &\sim  \int_0^t (t-r)^{2\alpha -2} ( \l (t-r)^{\alpha})^{-d} \prod_{j=1}^d x_j dr \nonumber\\
 &\sim \int_0^t (t-r)^{2\alpha -2 - \alpha d} dr,
\end{align}
which is finite if and only if $2\alpha-2 -\alpha d > -1$, i.e. if and only if $d=1$.\\

2) Assume $\alpha <1.$\\
Then we have 
\begin{align}
\Big(E_{\alpha,\alpha}(-\lambda t^{\alpha} |y|^2)\Big)^2 \sim (\lambda t^{\alpha})^{-2} |y|^{-4}
\end{align}
and by \eqref{11.38} we see that $\| G(u,\cdot) \|^2 = \infty$ in that case.\\
Combining this with Theorem \ref{th8.1} we get that the solution is not mild.\\

We summarise what we have proved in the following:
\begin{theorem}\label{mild2}
For all $H=(H_0,H_1, ... ,H_d) \in (\tfrac{1}{2},1)^{1+d}$ the following holds for the stochastic fractional Volterra equation \eqref{vol}:
\begin{itemize}
\item
(i)   If $\alpha = 1$ the solution $Y(t,x)$  is mild if and only if $d=1$.
\item
 (ii)   If $\alpha > 1$ the solution is mild if $d=1$ or $d=2$.
 \item
 (iii) If $\alpha < 1$ the solution is not mild for any $d$.
 \end{itemize}
\end{theorem}

\begin{remark}
    It is surprising that these results do not depend on the Hurst parameter $H=(H_0,H_1, ... ,H_d) \in (\tfrac{1}{2},1)^{1+d}$.
\end{remark}
\section{Application}
\subsection{Distillation columns}
Distillation process separates between liquid and vapor in many fields, where it works by applying  and removing a heat.
Heat causes components with lower boiling points and higher volatility to vaporize, leaving less volatile components in liquid form.

It separates binary mixtures. Many variables, such as column pressure, temperature, size, and diameter are determined by the properties of the feed and the desired products.  
\subsection{A model for the heat transfer in distillation columns}
To model the transfer of heat in such a mixture we use the following fractional heat equation for the temperature $T(t,x)$ at time $t$ and at the position $x$ in the $3-$dimensional space:

\begin{align} 
&\frac{\partial^{\alpha}}{\partial t^{\alpha}}T(t,x)=\lambda \Delta T(t,x)+\sigma T(t,x) \diamond W_{H}(t,x);\; (t,x)\in [0,\infty)\times [0,\infty)^3
\end{align}
   with boundary conditions
\begin{align}
    T(0,x)&=\delta_0(x)\text{ (the Dirac measure at  } x), \nonumber\\
    T(t,0)&=0 \text{ when } t >0.
  \end{align}  
By Theorem \ref{mult1} the solution of this fractional stochastic differential equation 
coincides with the solution of the following fractional Volterra equation: 
\begin{align} \label{vol2}
T(t,x)&=f(t,x) + \int_0^t \int_{0}^x  g(t-r,x-z) T(r,z) B_H(dr,dz),
\end{align}
where
\small{
\begin{align}
f(t,x)=(2\pi)^{-d} \int_{\mathbb{R}^3} e^{ixy} E_{\alpha}(- \l t^{\alpha} |y|^2) dy
=(2\pi)^{-3} \int_{\mathbb{R}^d} e^{ixy}\sum_{k=0}^{\infty} \frac{(- \l t^{\alpha} |y|^2)^k}{\Gamma(\alpha k +1)}dy, 
\end{align}
}
and
\small{
\begin{align}
    &g(t-r,x-z)= \sigma (2\pi)^{-3} (t-r)^{\alpha -1}\int_{\mathbb{R}^{3}}e^{i(x-z)y} E_{\alpha,\alpha}(-\l (t-r)^{\alpha}|y|^2) dy  
\end{align}
A solution formula for this Volterra equation is given in Theorem \ref{vol3}.\\

The question is:\\
\textit{For what values of $\alpha$ do we obtain vaporisation, and for what $\alpha$ do we obtain volatility? 
}\\

Fractional calculus have an important  role to describe real-world phenomena  called memory effect. The Caputo derivative is used because the memory effect and a constant function derivatives yield zero. Or the distillation need the memory effect so we use it to interpret the situation.
In our previous paper \cite{MO} we show that in the case where $\alpha< 1$ the travel time of the particles are longer than the standard diffusion, this case models also the removal of heat.
Where $\alpha>1$ the particles spread faster than in regular diffusion, in this case we obtain a vaporisation.

\section{Acknowledgments}
We are grateful to Yaozhong Hu for helpful comments.

\section{Appendix}
\subsection{Multiple-parameter  fractional Brownian motion}
Next we introduce the multiple-parameter \emph{fractional} Brownian motion and give some of its properties:
\begin{definition}\cite{OZ}
 We define n-parameter fractional Brownian motion $B_{H}(x)$; $x=(x_{1},x_{2},...,x_{d})\in \mathbb{R}^{n}$ with Hurst parameter $H=(H_{1},...,H_{n})\in (0,1)^{n}$  as a Gaussian process on $\mathbb{R}^{n}$ with mean
 \begin{equation*}
     \E[B_{H}(x)]=0\; \text{for all}\:\; x\in \mathbb{R}^{n}
 \end{equation*}
 and covariance 
 \begin{equation*} \label{cov}
     \E[B_{H}(x)B_{H}(y)]=\big(\frac{1}{2}\big)^{n}\prod_{j=1}^{n}\big( |x_{j}|^{2H_{j}}+|y_{j}|^{2H_{j}}-|x_{j}-y_{j}|^{2H_{j}}\big)
 \end{equation*}
 We also assume that 
 \begin{equation*}
     B_{H}(0)=0\;a.s. 
 \end{equation*}
\end{definition}
 From now, we assume that $\frac{1}{2}<H_{i}<1$ for $i=1,...,n$.\\
  \begin{definition}
For $x=(x_1,x_2, ... , x_n), y=(y_1,y_2, ... ,y_n) \in \R^n$ define 
\begin{align}
 \varphi(x,y)=\varphi_H(x,y)= \prod_{j=1}^{n}H_j(2H_j -1) |x_j - y_j |^{2H_j -2},
 \end{align}
 and let $L_{\varphi}^2(\R^n)$ be the set of all functions $f:\R^n \mapsto \R$ such that
 \begin{align}
 |f|_{\varphi}^2:=\int_{\R^n} \int_{\R^n} f(x) f(y) \varphi(x,y) dx dy < \infty.
\end{align}
Then $L_{\varphi}^2(\R^n)$ is a separable Hilbert space, with inner product
\begin{align}
    \<f,g\>:=\int_{\R^n} \int_{\R^n} f(x) g(y) \varphi(x,y) dx dy; \quad f,g \in L_{\varphi}^2(\R^n).
\end{align}
\end{definition}
The following \emph{fractional Ito isometry} is useful:
\begin{proposition} 
Let $f: \R^n \mapsto \R$ be a deterministic function. Then
\begin{align}\label{isometry}
    &E[(\int_{\R^n} f(x) dB_{H}(x))^2] = |f|_{\phi}^2 = \int_{\R^n} \int_{R^n} f(x)f(y) \varphi(x,y) dx dy=\nonumber\\
&=\prod_{j=1}^{n}H_j(2H_j -1) \int_{\R^n} \int_{R^n} f(x)f(y) |x_j - y_j |^{2H_j -2} dx dy ; \quad f \in L_{\varphi}^2(\R^n).
\end{align}
\end{proposition}

\subsection{Multiparameter fractional white noise}
We need to define the multiparameter fractional Hida test function space $\big(\S\big)_{H}$ and distribution space $\big(\S\big)^{*}_{H}$.
Let $\N_0=\{0,1,2,...\}$ denote the set of nonnegative integers, and let $\mathcal{J}$ denote the set of all (finite) multi-indices  
$\alpha = (\alpha_1, \alpha_2, ... , \alpha_m)$ with $\alpha_i \in \mathbb{N}_0; i=1,2, ... ,m; \quad m=1,2, ....$
\begin{definitionn}\cite{OZ}
Let $\mathcal{H}_{\alpha}$ be as in \cite {OZ},Section 2.
    \begin{itemize}
        \item (The multiparameter fractional Hida test function spaces)\\
        For $k\in \N$, define $(\S)_{H,k}$ to be the space of all expansions
        \begin{equation*}\label{psi}
            \psi(\omega)=\sum_{\alpha\in J}a_{\alpha}\mathcal{H}_{\alpha}(\omega)
        \end{equation*}
        such that 
        \begin{equation*}
           \| \psi \|^{2}_{H,k} :=\sum_{\alpha\in J}\alpha!a^{2}_{\alpha}(2\N)^{k\alpha}<\infty,
        \end{equation*}
        where $(2\N)^{\gamma}=\prod_{j}(2j)^{\gamma_{j}}$ if $\gamma=(\gamma_{1},...\gamma_{m})\in \mathcal{J}.$
        
        Define $(\S)_{H}=\bigcap_{k=1}^{\infty}(\S)_{H,k}$, equipped with the projective topology.
        \item (The multiparameter fractional Hida distribution spaces)\\
        For $q\in N$, let $(\S)^{*}_{H,-q}$ be the space of all formal expansions 
        \begin{equation*}\label{G}
            G(\omega)=\sum_{\beta\in J}b_{\beta}\mathcal{H}_{\alpha}(\omega),
        \end{equation*}
        such that 
        \begin{equation*}
            \|G\|^{2}_{H,-q}:=\sum_{\beta\in J}\beta!b_{\beta}^{2}(2N)^{-q\beta}<\infty.
        \end{equation*}
        Define $(\S)^{*}_{H}=\cup _{q=1}^{\infty}(\S)^{*}_{H,-q}$ with the inductive topology. Then $(\S)^{*}_{H}$ becomes the dual of $(\S)_{H}$ when the action of $G(\omega)=\sum_{\beta\in J}b_{\beta}\mathcal{H}_{\alpha}(\omega)\in (\S)^{*}_{H} $ on $\psi(\omega)=\sum_{\alpha\in J}a_{\alpha}\mathcal{H}_{\alpha}(\omega) \in (\S)_{H}$ is defined by
        $$\<\<G,\psi\>\>:=\sum_{\alpha\in J}\alpha!a_{\alpha}b_{\alpha}. $$
    \end{itemize}
\end{definitionn}
\begin{examplee}\cite{OZ}(Multiparameter fractional white noise)\\
Define for $y\in \R^{n}$
\begin{equation}
    W_{H}(y)=\sum_{i=1}^{\infty}\Bigg[\int_{\R^{n}}e_{i}(v)\phi(y,v)dv\Bigg]\mathcal{H}_{\epsilon(i)}(\omega).
\end{equation}
  Then as in \cite{HO},Example 3.6 we obtain that $W_{H}(y)\in (\S)^{*}_{H}$ for all $y$. Moreover $W_{H}(y)$ is integrable in $(\S)^{*}_{H}$ for $0\leq y_{i}\leq x_{i},\;\; i=1,...,d$ and 
  \begin{align}
      \int_{0}^{x_1}\int_{0}^{x_2} ... \int_{0}^{x_n}W_{H}(y)dy&= \sum_{i=1}^{\infty}\Bigg[\int_{0}^{x_1}\int_{0}^{x_2} ... \int_{0}^{x_n}\Big(\int_{\mathbb{\R}^{n}}e_{i}(v)\phi(y,v)dv\Big)dy\Bigg]\mathcal{H} _{\epsilon_{i}}(\omega)\nonumber\\
      &=B_{H}(x).
  \end{align}
  Therefore  $B_{H}(x)$ is differentiable with respect to $x$ in $(\S)^{*}_{H}$ and we have 
  \begin{equation} \label{whitenoise}
      \frac{\partial^{n}}{\partial x_{1},...,\partial x_{n}}B_{H}(x)=W_{H}(x)\;\; in\;(\S)^{*}_{H}.
  \end{equation}
  This justifies the name (multiparameter) fractional white noise for $W_{H}(x)$.
\end{examplee}
\begin{definitionn}\cite{OZ}
    Suppose $F(\omega)=\sum_{\alpha\in J}a_{\alpha}\mathcal{H}_{\alpha}(\omega)$ and $G(\omega)=\sum_{\beta\in J}b_{\beta}\mathcal{H}_{\beta}(\omega)$ both belong to $(\S)^{*}_{H}$. Then we define their Wick product $\Big(F \diamond G\Big)(\omega)$ by
    \begin{equation}
       \Big(F \diamond G\Big)(\omega)=\sum_{\alpha,\beta \in J}a_{\alpha}b_{\beta}\mathcal{H}_{\alpha+\beta}(\omega)=\sum_{\gamma\in J}\Bigg(\sum_{\alpha+\beta=\gamma}a_{\alpha}b_{\beta}\Bigg)\mathcal{H}_{\gamma}(\omega) 
    \end{equation}
\end{definitionn}
\begin{examplee}\cite{OZ}
    \begin{itemize}
        \item If $f,g \subset L^{2}_{\phi}(\R^{d})$ then
        \begin{equation}
            \Bigg(\int_{\R^{d}}f dB_{H}\Bigg)\diamond \Bigg(\int_{\R^{d}}gdB_{H}\Bigg)=\Bigg(\int_{\R^{d}}fdB_{H}\Bigg)\Bigg(\int_{\R^{d}}gdB_{H}\Bigg)-(f,g)_{\phi}
        \end{equation}
        \item If $f\in L^{2}_{\phi}(\R^{d})$ then\\
        $exp^{\diamond}(<\omega,f>):=\sum_{n=1}^{\infty} \frac{1}{n!}<\omega,f>^{\diamond n}$ converges in $(\S)^{*}_{H}$ and is given by 
        \begin{equation}
         \exp^{\diamond}\Big(<\omega,f>\Big)=\exp \Big(<\omega,f>-\frac{1}{2}|f|^{2}_{\phi}\Big).   
        \end{equation}
        
    \end{itemize}
\end{examplee}
\begin{definitionn}\cite{OZ}
    Suppose $y:\R^{d}\longrightarrow (\S)^{*}_{H}$ is a given function such that $x \mapsto y(x)\diamond W_{H}(x)$ is integrable in $(\S)^{*}_{H}$ for $x\in \R^{n}$.
    Then we define the multiparameter fractional stochastic integral (of Ito  type) of $y(x)$ by 
    \begin{equation}
        \int_{\R^{n}}y(x)dB_{H}(x)=\int_{\R^{n}}y(x)\diamond W_{H}(x)dx
    \end{equation}

    \end{definitionn}

\end{document}